\documentclass[numbers,webpdf,imaiai,authoryear]{ima-authoring-template}

\usepackage{subfig}
\usepackage{natbib}
\usepackage{float}
\usepackage{booktabs}
\usepackage{xcolor}

\graphicspath{{Fig/}}
\theoremstyle{thmstyleone}
\newtheorem{theorem}{\textbf{Theorem}}

\newtheorem{remark}{\textbf{Remark}}
\theoremstyle{thmstyletwo}%
\theoremstyle{thmstylethree}%

\numberwithin{equation}{section}

\begin{document}

\DOI{DOI HERE}
\copyrightyear{2025}
\vol{00}
\pubyear{2025}
\access{Advance Access Publication Date: Day Month Year}
\appnotes{Paper}
\copyrightstatement{Preprint submitted to arxiv.org}
\firstpage{1}

\title[Extremum Seeking for PDE Systems using Physics-Informed Neural Networks]{Extremum Seeking for PDE Systems using Physics-Informed Neural Networks}

\author{\textsc{Haojin Guo, Zongyi Guo*, Jianguo Guo}
\address{\orgdiv{Institute of Precision Guidance and Control}, \orgname{Northwestern Polytechnical University}, \orgaddress{\street{Xi'an}, \postcode{710072}, \state{Shaanxi}, \country{China}}}}

\author{\textsc{Tiago Roux Oliveira}
\address{\orgdiv{Dept. of Electronics and Telecommunication Engineering}, \orgname{State University of Rio de Janeiro}, \\ \orgaddress{\street{Rio de Janeiro}, \postcode{20550-900}, \state{RJ}, \country{Brazil}}}}

\corresp[*]{Corresponding author. Email:\href{email:guozongyi@nwpu.edu.cn}{guozongyi@nwpu.edu.cn}}

\authormark{H. Guo \textit{ET AL.}}

\received{Date}{0}{Year}
\revised{Date}{0}{Year}
\accepted{Date}{0}{Year}

\abstract{\textit{Abstract}: Extremum Seeking (ES) is an effective real-time optimization method for PDE systems in cascade with nonlinear quadratic maps. To address PDEs in the feedback loop, a boundary control law and a re-design of the additive probing signal are mandatory. The latter, commonly called “trajectory generation” or “motion planning,” involves designing perturbation signals that anticipate their propagation through PDEs. Specifically, this requires solving motion planning problems for systems governed by parabolic and hyperbolic PDEs. Physics-Informed Neural Networks (PINN) is a powerful tool for solving PDEs by embedding physical laws as constraints in the neural network’s loss function, enabling efficient solutions for high-dimensional, nonlinear, and complex problems. This paper proposes a novel construction integrating PINN and ES, automating the motion planning process for specific PDE systems and eliminating the need for case-by-case analytical derivations. The proposed strategy efficiently extracts perturbation signals, optimizing the PDE system.}


\maketitle

\section{Introduction}

Partial differential equations (PDEs) are fundamental in the study of dynamic systems, essentially for modeling complex behaviors. They describe how quantities evolve over time and space under various conditions. The Navier-Stokes equations in fluid mechanics, the heat conduction equation, and Maxwell’s equations in electromagnetism are classical examples, widely applied in fields like aerospace, building insulation, wireless communication, and electromagnetic wave propagation \citep{2,4,5,34}.

Extremum Seeking (ES) is a widely used real-time optimization method for nonlinear dynamic systems, designed to identify the optimal objective function by adjusting the input parameters \citep{10,13,31}. In nonlinear dynamic systems, ES typically uses high-frequency small-amplitude perturbation signals and low-pass filters to estimate the gradient of the objective function, guiding the system toward the extremum \citep{8,33}. Adaptive optimization techniques have rapidly advanced, and in the field of control, \citep{60} implemented an adaptive controller for spacecraft attitude control, while \citep{61} designed an adaptive fault-tolerant controller to address spacecraft failures. In the chemical engineering domain, \citep{62} applied adaptive optimization to the thermal performance parameters of spacecraft. Similarly, ES has become a key optimization tool with applications to areas like spacecraft attitude control and chemical process optimization \citep{35}. As noted by Liu et al. \citep{11}, the global extremum seeking control method was introduced to improve both transient and steady-state responses. Furthermore, stochastic perturbation and adaptive perturbation-demodulation loops have proven effective in ensuring convergence in the presence of system disturbances and nonstationary noise \citep{12,15}. In this context, Oliveira et al. \citep{14} extended ES to infinite-dimensional systems, applying it to transport hyperbolic and diffusion parabolic PDEs with delays. Silva et al. \citep{16} proposed a stochastic ES method for multivariable systems with delays, incorporating predictor feedback to ensure stability and convergence.

On the other hand, traditional numerical methods are widely used to solve PDEs, but they face significant limitations when dealing with complex systems \citep{3}. Methods such as finite difference, finite element, and spectral methods rely on discretization and require specific adjustments to ensure accuracy \citep{6}. However, when handling high-dimensional, nonlinear, or complex boundary condition systems, these methods often exhibit inefficiency, especially during long-term simulations of dynamic systems \citep{7,22}. As a result, there is a pressing need for a universal solution
strategy capable of handling multiple PDE systems \citep{32}.

PINNs have been widely used in solving PDEs, yielding significant results by embedding physical laws into the neural network training, thus improving accuracy \citep{23}. Recent studies have focused on optimizing the PINN framework, such as using transfer learning to address multi-fidelity data challenges \citep{21}, composite neural networks for inverse PDE solutions \citep{20}, and Euler-PINNs integrated with PID controllers for dynamic systems under uncertainty \citep{63}. Additionally, several studies have introduced random and adaptive perturbation loops to handle uncertainty and noise in PDEs \citep{15,5}. However, existing methods still face limitations in solving disturbance signals across different PDE systems. To address this, we propose a unified PINN-based model that simplifies probing/perturbation signal design by adjusting inputs and loss functions, offering a new approach for ES in PDE systems. The key idea is to use a single neural network framework that efficiently solves multiple PDE systems and extracts the additive probing signal usually required in ES, simplifying the complexity of analytically obtaining such signals in the current metodology of ES for PDE systems \citep{36}. Our contributions can be summarized as follows:

\begin{itemize}
    
\item To summarize in a tutorial perspective a small portion of the vast space of possibilities of designing and analysing ES algorithms for infinite-dimensional systems governed by  diffusion parabolic PDEs. It is the exhaustive consideration of this baseline problem that make possible the various extensions—multivariable maps \citep{37}, dynamical plants \citep{38}, nonconstant \citep{39} and uncertain delays \citep{40}, distributed delays \citep{41}, reaction–advection–diffusion equations \citep{42} or wave models \citep{43}, cascades of PDEs \citep{13}, etc.—of the use of delay- and PDE-compensated ES algorithms, some of them discussed at the end of the article.
\item The unified solving framework based on PINN proposed in this study provides a single solving method for different types of PDE systems. Through this framework, there is no need to design separate solving methods for each PDE system, which significantly improves solving efficiency.
\item For the ES strategies of different PDE systems, this approach does not require changing the solving strategy for obtaining the probing/perturbation signals, thus avoiding the issue of frequently switching solving methods in traditional approaches.

\end{itemize}

Taken all together, this paper presents as the main novelty a method for solving complex PDE systems using a unified PINN framework, particularly suitable for optimization problems involving probing/perturbation signals, highlighting the potential of the PINN method in real-time optimization. Although, one could explore how PINNs or other machine-learning techniques could automate this process for a broader class of PDEs, we have eliminated the need for case-by-case analytical derivations for PDE motion planning and its application in ES as the main benefit. We believe that integrating ES control with machine learning offers a promising and innovative approach, enabling the experimental discovery of equivalent solutions through PINNs. This fusion could significantly extend the applicability of ES to even more complex dynamical systems while preserving rigorous convergence guarantees.

The paper is organized as follows. Section 2 presents the problem statement, introducing the basic ES model and the approach for handling particularly diffusion-dominated parabolic PDE systems. Section 3 is devoted to the methodology, describing the combination of PINN framework for trajectory generation and the boundary feedback ES design via backstepping transformation for PDE compensation. Section 4 provides a stability proof, ensuring the necessary step-by-step process to properly sequence them, with only one correct sequence, and to not fail in producing a valid proof to the considered case as well as across various other problems involving differents types of PDEs. The design principles of the ES compensator in Section 3 and the stability proof approach discussed in Section 4 are primarily based on the methods presented in \citep{14,36}. Section 5 presents numerical simulations, illustrating the method’s performance in practical applications. Finally, Section 6 summarizes the contributions and discusses future research directions. 

\section{Problem statement}

For control systems consisting of PDEs, we consider a real-time optimal control scheme, which aims to adjust the input to find and maintain the optimal value $\Theta^*\in \mathbb{R}$ of an unknown nonlinear static mapping $Q:\mathbb{R}\rightarrow \mathbb{R}$. The optimal input of this mapping is unknown and denoted by $\Theta\in \mathbb{R}$, while the output is measurable as $y\in \mathbb{R}$. To simplify the problem, we consider the maximization case. 

\subsection{Actuation dynamics and demodulation signals}

We take distributed PDE within the class of parabolic PDE as an example, and define the actuator as \( \theta(t) \in \mathbb{R} \) and the propagated actuator as \( \Theta(t) \in \mathbb{R} \). The state of the field $\alpha(x,t)$ satisfies as follows:
\begin{align}
\label{01}
&\Theta(t)=\alpha(0,t) \\
\label{eq:heat_eqn_1_CP8}
&\partial_t\alpha(x,t)=\partial_{xx}\alpha(x,t),\quad x\in(0,D) \\
&\partial_x\alpha(0,t) = 0 \\
\label{04}
&\alpha(D,t)= \theta(t),
\end{align}
where $\alpha: [0,D]\times \mathbb{R}_+ \rightarrow \mathbb{R}$ and $D>0$ is the diffusion coefficient of the diffusion PDE.
The measurement is defined by the unknown static map with input \eqref{01}, For simplicity, we assume the unknown nonlinear static map is quadratic, \textit{i.e.}, such that 
\begin{align} \label{05}
	y(t) = Q(\Theta(t)), \quad y(t) = y^* + \frac{H}{2}(\Theta(t) - \Theta^*)^2.
\end{align} 

We consider the unknown optimal input \( \Theta^* \in \mathbb{R} \) and output \( y^* \in \mathbb{R} \), with the unknown Hessian \( H < 0 \). Following the approach in \citep{25}, a perturbation signal \( S(t) \) with amplitude \( a \) and frequency \( \omega \) is introduced to estimate the gradient of the objective function. Additionally, we define the state parameter \( \beta(x,t) \), where \( \beta : [0,D] \times \mathbb{R}_+ \to \mathbb{R} \).
\begin{align}
\label{06}
&S(t):=\beta(D,t) \\ 
&\partial_t\beta(x,t)=\partial_{xx}\beta(x,t), \quad x\in(0,D) \\
&\partial_x\beta(0,t)=0 \\
\label{09}
&\beta(0,t)=a\sin(\omega t).
\end{align}

Traditional numerical methods discretize the PDE over the entire spatial domain with spatial nodes $N_x$. To improve accuracy, a larger $N_x$ and smaller time step $\Delta t$ are required, leading to a computational workload proportional to both $N_x$ and $N_t$, increasing overall complexity.
\begin{equation} \label{10}
    T(N_x \times N_t) = O(N_x \times N_t).
\end{equation}
Therefore, for different PDE systems, extremum seeking faces high computational costs and low flexibility in obtaining the perturbation signal $S(t)$, making it difficult to meet the needs of fast optimization \citep{17}. To address this, we propose using a unified PINN framework to solve for the perturbation signal.

\subsection{Estimation Errors and PDE-Error Dynamics}

To detect and demodulate the signal, and drive the unknown controller $\theta(t) \equiv \Theta^*$ to the optimal $y^*$, we estimate $\Theta^*$, with the estimation error defined between the optimal and propagated input variables:
\begin{align} \label{11}
&\hat \theta(t) = \theta(t) - S(t), \quad \hat \Theta(t) = \Theta(t) - a \sin (\omega t), \\
\label{12}
&\tilde{\theta}(t) := \hat{\theta}(t)-\Theta^*, \quad \vartheta(t) := \hat{\Theta}(t)-\Theta^*.
\end{align}
The relationship between the propagation estimation error $\vartheta(t)$, the input $\Theta(t)$, and the static mapping optimizer $\Theta^*$ is derived from Equations \eqref{11} and \eqref{12}:
\begin{equation} \label{13}
\vartheta(t)+a\sin(\omega t) = \Theta(t)-\Theta^*.
\end{equation}
Let us define
\begin{equation} \label{14}
\bar{\alpha}(x,t) := \alpha(x,t) - \beta(x,t) - \Theta^*, \quad \dot{\tilde{\theta}} = U(t),
\end{equation}
where $u(x,t) := \partial_t \bar{\alpha}(x,t) = \partial_{xx} \bar{\alpha}(x,t)$. 
Using Equations \eqref{01}-\eqref{04} and \eqref{11}-\eqref{14}, combined with Equations \eqref{06} and \eqref{09}, and assume that $\bar{\alpha}, \bar{\alpha}_t, \bar{\alpha}_x, \bar{\alpha}_{xx} \in \mathcal{C}^1([0,D] \times \mathbb{R}_+)$, see Fig.~\ref{ESloop}, the propagated error dynamics can then be written as:
\begin{align}
\label{15}
&\dot{\vartheta}(t)=u(0,t), \\ 
&\partial_t u(x,t)=\partial_{xx}u(x,t), \quad x\in(0,D) \\
&\partial_x u(0,t)=0, \\
&u(D,t)=U(t).
\label{18}
\end{align}
It is essential to define the demodulation signal \( N(t) \), which is employed to estimate the Hessian of the static map by multiplying it with the output \( y(t) \), as detailed in \citep{26}. Meanwhile, the dither signal \( M(t) \), which estimates the gradient, remains unchanged from the basic ES approach:
\begin{figure*}[htbp] 
    \centering
    \includegraphics[width=0.80\linewidth]{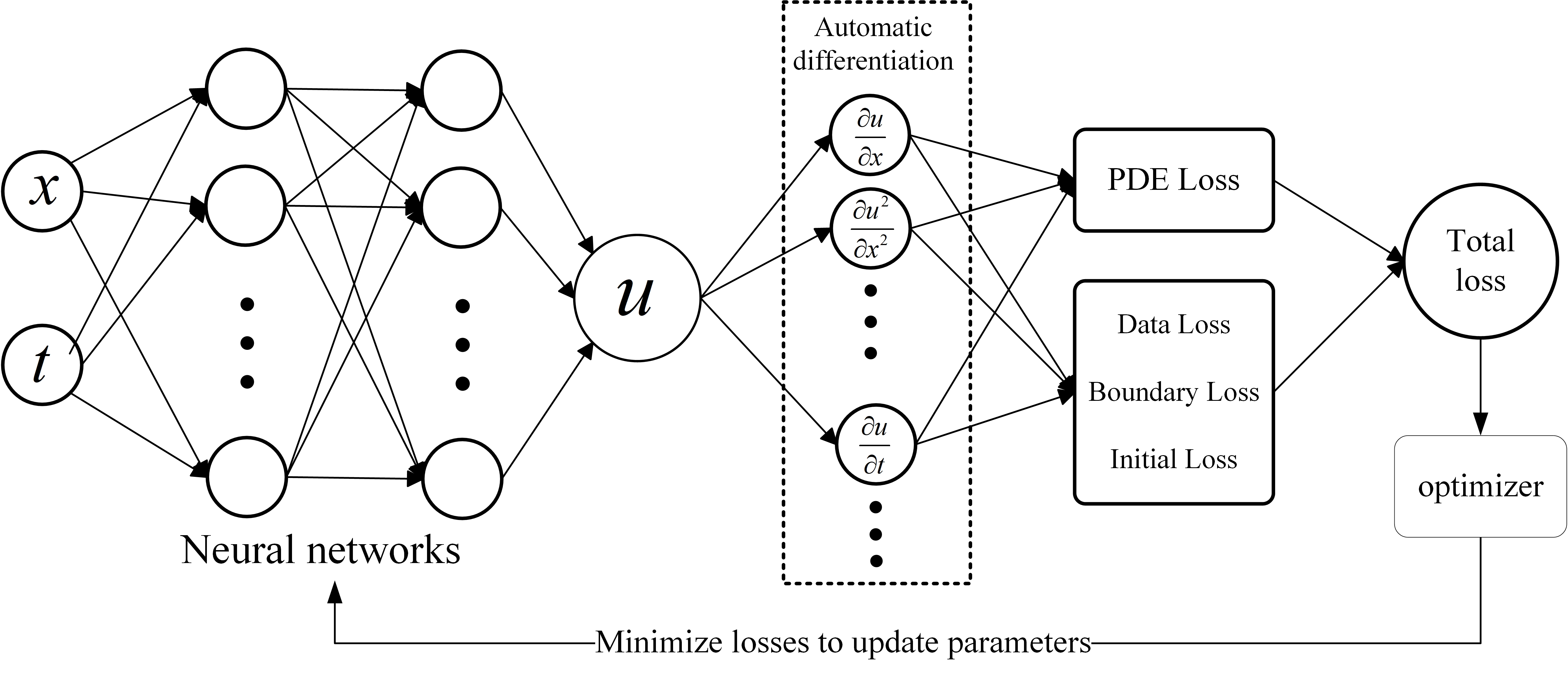} 
    \caption{Schematic Diagram of the PINN Architecture}
    \label{fig:horizontal-image}
\end{figure*}
\begin{align}
\label{19}
&\hat{H}(t) = N(t)y(t)\quad \text{with} \quad N(t) = -\frac{8}{a^2}\cos(2\omega t),\\
\label{20}
&G(t) = M(t)y(t)\quad \text{with} \quad M(t) = \frac{2}{a}\sin(\omega t).
\end{align}

\section{Methodology}

\subsection{Physics-Informed Neural Networks (PINN)}

The principle of PINN is to approximate the solution $u(x,t)$ using a neural network $u_{\theta}(x,t)$, where $\theta$ represents the network's trainable parameters. The network takes $x$ and $t$ as inputs and outputs the corresponding solution approximation.
\begin{equation} \label{21}
    u_\theta : (x,t) \mapsto u_\theta(x,t) \approx u(x,t).
\end{equation}
The latent solution $u(x,t)$ is approximated by a deep neural network $u_{\theta}(x,t)$, which allows direct evaluation of the PDE residuals in the loss function:
\begin{equation} \label{22}
    L_{\text{PDE}}(\theta) = \frac{1}{N} \sum_{i=1}^{N} \left[ \mathcal{L}[u_{\theta}(x_i, t_i)] + \mathcal{N}[u_{\theta}(x_i, t_i)] \right]^2.
\end{equation}
The PDE residual represents the difference between the neural network's predicted solution $u_{\theta}(x,t)$ and the physical laws. The PINN loss function $L(\theta)$ consists of the PDE error $L_{\text{PDE}}(\theta)$, data error $L_{\text{data}}(\theta)$, initial condition error $L_{\text{IC}}(\theta)$, and boundary condition error $L_{\text{BC}}(\theta)$, and is given by:
\begin{align} 
\label{23}
&L(\theta) = L_{\text{data}}(\theta) + L_{\text{PDE}}(\theta) + L_{\text{BC}}(\theta) + L_{\text{IC}}(\theta)\\
&L_{\text{data}} = \frac{1}{N_u} \sum_{i=1}^{N_u} \left| u_{\theta}(x_u^i, t_u^i) - u_{\text{data}}^i \right|^2\\
\label{25}
&L_{\text{IC}} = \frac{1}{N_{\text{IC}}} \sum_{i=1}^{N_{\text{IC}}} \left| u_{\theta}(x_{\text{IC}}^i, t_{\text{IC}}^i) - h(x_{\text{IC}}^i) \right|^2\\
\end{align}
\begin{align} 
\label{26}
&L_{\text{BC}} = \frac{1}{N_{\text{BC}}} \sum_{i=1}^{N_{\text{BC}}} \left| u_{\theta}(x_{\text{BC}}^i, t_{\text{BC}}^i) - g(x_{\text{BC}}^i, t_{\text{BC}}^i) \right|^2.
\end{align}
When dealing with different PDE systems, we utilize the PINN framework to obtain the perturbation signal \( S(t) \), with the process following Equations \eqref{06} to \eqref{09}:
\begin{equation} \label{27}
a\sin(\omega t) \xrightarrow{\text{PINN}} S(t)
\end{equation}

By embedding physical laws into the loss function, PINNs eliminate the need for grid discretization, making them applicable to a wide range of PDEs, which simplifies the solution process and improves computational efficiency, whether the PDE is linear, nonlinear, or has complex boundary conditions.

\subsection{Control Design and Closed-Loop Error Dynamics}
We adopt the ES strategy depicted in Fig~\ref{ESloop} and design a compensator accordingly. We consider the PDE-ODE cascade \eqref{15}-\eqref{18}, and, as in \citep{24}, apply the backstepping transformation.
\begin{align} 
\label{28}
w(x,t) = u(x,t)-\int_{0}^{x}q(x,r)u(r,t)dr-\gamma(x)\vartheta(t),
\end{align}
with gain kernels\index{kernel PDE} $q(x,r) = \bar{K}(x-r)$ and $\gamma(x) = \bar{K}$, which transforms \eqref{23}-\eqref{26} into the target system
\begin{align} 
\label{29}
\dot{\vartheta}(t) &= \bar{K}\vartheta(t) + w(0,t) \\ 
\label{eq:target_2}
\partial_t w(x,t) &= \partial_{xx}w(x,t), \quad x \in (0,D) \\ 
\label{eq:target_3}
\partial_x w(0,t) &= 0 \\ 
\label{32}
w(D,t) &= 0,
\end{align}

with $\bar{K}< 0$.

As the target system \eqref{29}-\eqref{32} is exponentially stable, the controller that compensates the diffusion process can be obtained by evaluating the backstepping transformation \eqref{28} at $x=D$:
\begin{align}  
\label{33}
U(t) = \bar{K}\vartheta(t)+\bar{K}\int_{0}^{D}(D-r)u(r,t)dr.
\end{align}
Since $\vartheta(t)$ cannot be measured directly, We extend the method to the parabolic PDE case by calculating the average version of the gradient in \eqref{20} and the Hessian estimate in \eqref{19} to resolve this challenge. 

The compensation system that processes the signal \( G_\text{av}(t) \) for the parabolic PDE system’s ES controller as follows:
\begin{align} 
\label{34}
G_\text{av}(t) = H\vartheta_\text{av}(t).
\end{align}
If a quadratic map as in \eqref{05} is considered, and according to \eqref{34}, we average \eqref{33} and choose \(\bar{K} = KH\) with \(K > 0\) and the introduced unknown Hessian \(H\) of the static map, Then, by computing the averaged gradient and Hessian estimates, we obtain the desired result:
\begin{figure*}[htbp] 
\centering
\includegraphics[width=1.0\linewidth]{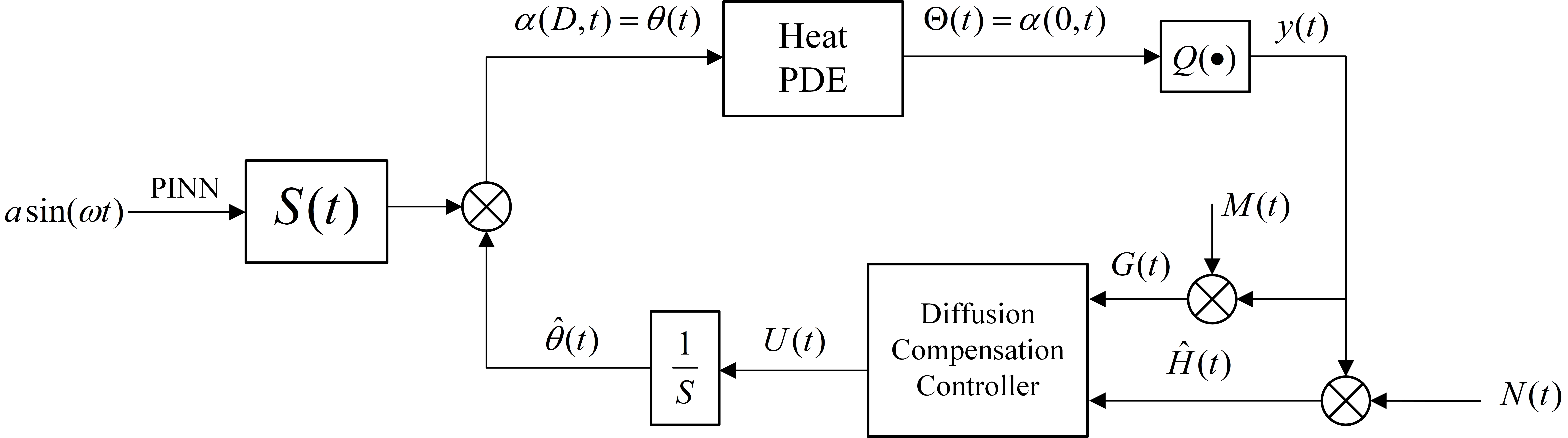} 
\caption{Schematic Diagram of the ES Architecture}
\label{ESloop}
\end{figure*}
\begin{align}
\label{35}
U_{\text{av}}(t) &= KG_{\text{av}}(t) + K\hat{H}_\text{av} \int_{0}^{D} (D-r) u_{\text{av}}(r,t) \, dr.
\end{align}
    
A low-pass filter is introduced into the controller to rewrite the error dynamics as an infinite-dimensional system, with $U(t)$ included as a state variable. Consequently, an average-based infinite-dimensional control law is derived to compensate for the diffusion process
\begin{align} 
\label{36}
U(t) =\mathcal{T}\left\lbrace K \left[G(t) + \hat{H}(t)\int_{0}^{D} (D-r)u(r,t)dr\right]\right\rbrace,
\end{align}
with the low-pass filter operator
\begin{align} 
\label{37}
\mathcal{T}\{\varphi(t)\}=\mathcal{L}^{-1}\left\{\frac{c}{s+c} \right\}*\varphi(t),\quad \varphi(t):\mathbb{R}_+\rightarrow \mathbb{R},
\end{align} 
where the corner frequency \( c > 0 \) is chosen later, \( \mathcal{L}^{-1}\{\cdot\} \) denotes the inverse Laplace transform, and \( * \) represents the convolution operator. Note that the feedback law \eqref{36} depends on the full PDE state information \( u(x,t) \), where \( x \in [0,D] \). Since \( u(x,t) \) is an artificial state used for analysis in the subsequent discussion, the feedback law \eqref{36} should be described in terms of the plant state \( \alpha(x,t) \) and dither signals. Applying integration by parts to the integral of \( \partial_t \beta(x,t) = \partial_{xx} \beta(x,t) \), as associated with \eqref{06}-\eqref{09}, the feedback law \eqref{36} is rewritten as follows:
\begin{equation} 
\label{39}
U(t) = \mathcal{T}\left\lbrace K \left[G(t) + \hat{H}(t) \left( \int_{0}^{D} (D-r) \partial_{t} \alpha (r,t)dr - (S(t) - a \sin (\omega t)) \right)\right]\right\rbrace.
\end{equation}
In practice, it is more reasonable to measure \( \Theta(t) \), the input of the static map, than to measure \( \partial_t \alpha(x,t) \), which is the time-derivative of the spatially distributed plant state. In this case, the control law \eqref{39} can be rewritten with the help of the diffusion equation \( \partial_t \alpha(x,t) = \partial_{xx} \alpha(x,t) \) and integration by parts, as applied to \eqref{01}-\eqref{04} and \eqref{11}-\eqref{12}:

\begin{align} \label{40}
U(t) =\mathcal{T}\left\lbrace K \left[G(t) + \hat{H}(t)\left(\hat{\theta}(t)-\Theta(t) + a \sin \left(\omega t \right) \right)\right]\right\rbrace.
\end{align}

We need to compute the error dynamics to assist in the stability analysis. With the approximated gradient and Hessian depending on $\vartheta(t)$ and $H$, the infinite-dimensional control law \eqref{36} is given by
\begin{align} \label{41}
U(t) =\ & \mathcal{T}\left\lbrace K\left[y^* + \frac{H}{2}\vartheta^2(t) + Ha\sin(\omega t)\vartheta(t) \right.\right. \nonumber\\
& + \left.\left. \frac{Ha^2}{2}\sin^2(\omega t)\right] \times \left[\frac{2}{a}\sin(\omega t) - \frac{8}{a^2}\cos(2\omega t)\int_{0}^{D}(D-r)u(r,t)dr\right]\right\rbrace. 
\end{align}
Substituting~(\ref{41}) into~(\ref{18}), we can write the error-dynamics \eqref{15}-\eqref{18} as 
\begin{align}
    \label{42}
    &\dot{\vartheta}(t) = u(0,t), \\
    &\partial_t u(x,t) = \partial_{xx} u(x,t), \quad x \in (0, D), \\
    &\partial_x u(0,t) = 0, \\
    &u(D,t) = \mathcal{T}\left\lbrace K \left[ Hf_1(t) + H\vartheta(t) + \sin(\omega t) f_2(t) \right.\right. \nonumber \\
    &\qquad \quad - \cos(2\omega t) f_3(t) - \sin(3\omega t) f_4(t) + \cos(4\omega t) f_5(t) \left. \left. \right] \right\rbrace,
    \label{45}
    \end{align}  
where 
\begin{eqnarray}
f_1(t) &=& \int_{0}^{D}(D-r)u(r,t)dr,\\ f_2(t)&=&(2y^*+H\vartheta^2(t)+4\vartheta(t)f_1(t)+3a^2H/4)/a,\\
f_3(t) &=& 8y^*f_1(t)/a^2 + H\vartheta(t)+2H + 4H\vartheta^2(t)f_1(t) / a^2,\\
f_4(t) &=& 4H\vartheta(t)f_1(t)/a+Ha/4,\\
f_5(t)&=& Hf_1(t).
\end{eqnarray}
\begin{remark} \label{remark1}
For the hyperbolic PDE system, such as transport PDEs and unlike the parabolic PDE system, there is a known delay \( D_t > 0 \) in the actuation or measurement process, which introduces differences in the design process. Specifically, this delay results in a static map for \( y^\prime(t) \) and demodulation signals \( \beta^\prime(0,t) \), where \( y^\prime(t) \) is given by the expression \( y^\prime(t) = y^* + \frac{H}{2}(\Theta(t-D_t) - \Theta^*)^2, \) and \( \beta^\prime(0,t) \) is represented as \( \beta^\prime(0,t) = a\sin(\omega (t+D_t)). \) In this case, due to the presence of the delay, the error dynamics for hyperbolic PDE systems are described by the equation \( \dot{\tilde{\theta}}^\prime(t - D_t) = U(t - D_t). \) It is essential for the compensation system that processes the signal \( G_\text{av}^\prime(t) \) for the hyperbolic PDE system’s ES controller to be designed as \( G_\text{av}^\prime(t) = H\vartheta_\text{av}(t - D_t) \) by means the introduction of predictor feedback \citep{ref}.
\end{remark}

\section{Stability and convergence analysis}

The following theorem summarizes the stability and convergence properties of the error-dynamics \eqref{42}-\eqref{45} and is proven in this section.
\begin{theorem} \label{theo}
    For a sufficiently large $c>0$, there exists some $\bar{\omega}(c)>0$, such that $\forall\omega > \bar{\omega}$, the error-dynamics \eqref{42}-\eqref{45} with states $\vartheta(t)$, $u(x,t)$, has an unique exponentially stable periodic solution in $t$ of period $\Pi:=2\pi/\omega$, denoted by $\vartheta^\Pi(t),u^\Pi(x,t)$, satisfying $\forall t\ge 0$:
\end{theorem}
\begin{equation} \label{401}
\left(\left|\vartheta^\Pi(t)\right|^2+\left\|u^\Pi(t)\right\|^2+\left\|\partial_x u^\Pi(t)\right\|^2 + \left|u^\Pi(D,t)\right|^2\right)^{1/2}\le \mathcal{O}\left(1/\omega\right).
\end{equation}
Furthermore,
\begin{align}
\label{402}
&\limsup\limits_{t\rightarrow \infty} |\theta (t)-\Theta^*| = \mathcal{O}\left(|a|e^{D\sqrt{\omega/2}}+1/\omega\right), \\
\label{403}
&\limsup\limits_{t\rightarrow \infty}|\Theta(t)-\Theta^*| = \mathcal{O}\left(|a|+1/\omega\right), \\
\label{404}
&\limsup\limits_{t\rightarrow \infty}|y(t)-y^*| = \mathcal{O}\left(|a|^2+1/\omega^2\right).
\end{align}
\textbf{\textit{Proof.}}~The proof is organized into Steps 1–6, following a similar structure to the approach outlined in \citep{14}.

\subsubsection*{\textbf{Step 1:}\ \ \textit{ Average Closed-Loop System}}
The average version of the system \eqref{42}-\eqref{45} for $\omega$ large is
\begin{align}
\label{405}
&\dot{\vartheta}_{\text{av}}(t)=u_{\text{av}}(0,t), \\
\label{406}
&\partial_t u_{\text{av}}(x,t)=\partial_{xx} u_{\text{av}}(x,t),\quad x\in[0,\ D] \\
\label{407}
&\partial_x u_{\text{av}}(0,t)=0, \\
&\frac{d}{dt}u_{\text{av}}(D,t)=-cu_{\text{av}}(D,t)
\!+\!cKH\left[\vartheta_{\text{av}}(t)\!+\!\int_{0}^{D}(D-r)\!u_{\text{av}}(r,t)dr\right].
\label{408}
\end{align}
The low-pass filter is represented in state-space form. To derive Equations \eqref{42}-\eqref{45}, the terms in Equation \eqref{45} that depend on the sine or cosine of the argument $k\omega$ for $k = 1, \ldots, 4$ are set to zero, as demonstrated in the derivation of Equation \eqref{34}.

\subsubsection*{\textbf{Step 2:}\ \ \textit{Backstepping Transformation into Target System}}

The backstepping transformation, as detailed in \citep{24}, is applied in \eqref{28} to facilitate the solution process.

\begin{equation}
w(x,t) = u_{\text{av}}(x,t) - KH\left[\vartheta_{\text{av}}(t) + \int_{0}^{x}(x - r) u_{\text{av}}(r,t) \, dr \right]
\label{409}
\end{equation}
maps the average error-dynamics \eqref{405}-\eqref{408} into the exponentially stable target system
\begin{align}
    \label{410}
    &\dot{\vartheta}(t) = u(0,t), \\
    &\partial_t u(x,t) = \partial_{xx} u(x,t), \quad x \in (0, D), \\
    &\partial_x u(0,t) = 0, \\
    &u(D,t) = \mathcal{T}\left\lbrace K \left[ H f_1(t) + H \vartheta(t) + \sin(\omega t) f_2(t) \right. \right. \nonumber \\
    &\qquad \quad - \cos(2\omega t) f_3(t) - \sin(3\omega t) f_4(t) + \cos(4\omega t) f_5(t) \left. \left. \right] \right\rbrace.
    \label{413}
\end{align}
Equations \eqref{410}-\eqref{413} are obtained by substituting the inverse backstepping transformation into the system, leading to the derived relations:

\begin{equation}
u_{\text{av}}(x,t) = w(x,t) + K H e^{K H x} \vartheta_{\text{av}}(t) + K H \int_{0}^{x} \left( e^{K H (x - r)} - 1 \right) w(r,t) \, dr.
\label{414}
\end{equation}
We first introduce the average error-dynamics as given by \eqref{405}-\eqref{408}. With the boundary Equation \eqref{410} and the backstepping transformation \eqref{409}, we obtain the relationship
\begin{align}
w(D,t) = -\frac{1}{c}\partial_t u_\text{av}(D,t), \quad \partial_t w(D,t) = \partial_t u_{\text{av}}(D,t) - K Hu_{\text{av}}(D,t).
\label{415}
\end{align}
Combined with \eqref{414} and \eqref{415}, leads to the result in \eqref{413}.

\subsubsection*{\textbf{Step 3:}\ \ \textit{Exponential Stability of the Target System}} 
Consider the Lyapunov-Krasovskii functional, we have
\begin{align}\label{417}
\Upsilon(t) = \frac{\vartheta_{\text{av}}^2(t)}{2} + \frac{a}{2} \| w(t) \|^2 + \frac{b}{2} \| \partial_x w(t) \|^2 + \frac{d}{2} w^2(D,t).
\end{align}
where \( a, b, d > 0 \). Define \( \lambda := -KH \) with \( \lambda > 0 \). Taking the time derivative of \eqref{417} and using integration by parts, we obtain the following result for the solution of the target system \eqref{410}-\eqref{413}.
\begin{align} \label{418}
\dot{\Upsilon}(t) = &\ -\lambda \vartheta_{\text{av}}^2(t) + \vartheta_{\text{av}}(t) w(0,t) + a w(D,t) \partial_x w(D,t) \nonumber \\
    &\ - a \|\partial_x w(t)\|^2 + b \partial_x w(D,t) \partial_t w(D,t) - b \|\partial_{xx} w(t)\|^2 \nonumber \\
    &\ + d w(D,t) \partial_t w(D,t).
\end{align}
By applying Young's, Poincaré's, Agmon's, and Cauchy-Schwarz inequalities to \eqref{418}, and selecting \( a = \frac{c - \lambda}{8D\lambda^3} \), \( b = \frac{1}{8D\lambda^3} \), and \( d = 1 \), we obtain the following result.
\begin{align} \label{419}
\dot{\Upsilon}(t) \le &\ -\frac{\lambda}{4} \vartheta_{\text{av}}^2(t) + (c_1^* - c) w^2(D,t) + (c_2^* - c) \|\partial_x w(t)\|^2 - \frac{1}{32 D \lambda^3} \|\partial_{xx} w(t)\|^2,
\end{align}
with 
\begin{align} \label{420}
&c_1^* = \frac{3}{2}\lambda^3 + \lambda + \frac{1+2D}{\lambda}+2D\lambda \zeta(D), \quad c_2^* = \lambda + 8D\lambda^3 \left[\frac{4D^2+1}{\lambda}+4D^2\lambda \zeta(D)\right],
\end{align}
and \( \zeta(D) = \int_{0}^{D} \left( e^{-\lambda (D - r)^2} - 1 \right) \, dr \). Therefore, from \eqref{419}, if \( c \) is chosen such that \( c > \max\left\lbrace c_1^*, c_2^* \right\rbrace \), we obtain for some \( \mu > 0 \):
\begin{align} \label{421}
\dot{\Upsilon}(t) \le - \mu \Upsilon(t).
\end{align}
Thus, the target system (4.12)-(4.13) is exponentially stable in the \( H_1 \)-norm
\begin{align} \label{422}
&\left(\vartheta^2_{\text{av}}(t)+\|w(t)\|^2+\|\partial_x w(t)\|^2+ w^2(D,t)\right)^{1/2},
\end{align}
\textit{i.e.}, in the transformed variable $(\vartheta_{\text{av}}, w)$.

\subsubsection*{\textbf{Step 4:}\ \ \textit{Exponential Stability Estimate (in $\mathcal{H}_1$-norm) of the Average Error-Dynamics}}
We define 
\begin{align} \label{423}
	\Psi(t) =&\ \vartheta^2_{\text{av}}(t)+\|u_{\text{av}}(t)\|^2 +\|\partial_x u_{\text{av}}(t)\|^2 + u^2_{\text{av}}(D,t).
\end{align}
Then, there are bounds for the Lyapunov-Krasovskii functional \eqref{417} with respect to \( \Psi(t) \).

\begin{equation} \label{424}
{\rho}\Psi(t) \le \Upsilon(t) \le \bar{\rho}\Psi(t), \quad {\rho} = {\tau}{\sigma},\quad \bar{\rho} = \bar{\tau}\bar{\sigma}, 
\end{equation}
with ${\sigma} = \min\left \lbrace \frac{1}{2},\ \frac{a}{2},\ \frac{b}{2}\right\rbrace, \bar{\sigma} = \max\left \lbrace \frac{1}{2},\ \frac{a}{2},\ \frac{b}{2}\right\rbrace$ and suitable \( \bar{\tau}, \tau \). From \eqref{423} and the exponential stability of the target system \eqref{420}, it follows that
\begin{align} \label{426}
	\Psi(t) \le \frac{\bar{\rho}}{{\rho}} e^{-\mu t} \Psi(0).
\end{align}
This completes the proof of the exponential stability of the average error dynamics \eqref{405}-\eqref{408} in the sense of the \( \mathcal{H}_1 \)-norm, with respect to \( \Psi^{1/2}(t) \), in the variables \( (\vartheta_\text{av}, u_\text{av}) \).

\subsubsection*{\textbf{Step 5:}\ \ \textit{Invoking the Averaging Theorem for Infinite-Dimensional Systems}}
To satisfy the conditions of the averaging theorem for infinite-dimensional systems, the error-dynamics \eqref{42}-\eqref{45} has to be in the form
\begin{equation} \label{427}
\dot{z}(t) = \Gamma z(t)+J(\omega t,z(t)).
\end{equation}
We utilize \( \Gamma \) and \( J(\omega t, z(t)) \) as defined in the theorem, and, through the state transformation in \eqref{42}-\eqref{45} with \( v(x,t) = u(x,t) - U(t) \), derive the error dynamics under homogeneous boundary conditions
\begin{align}
    \label{428}
    &\dot{\vartheta}(t) = v(0,t) + U(t), \\
    \label{429}
    &\partial_t v(x,t) = \partial_{xx} v(x,t) - \phi(\vartheta,v,U,t), \quad x \in (0, D) \\    
    \label{429_2}
    &\partial_x v(0,t) = 0, \\
    \label{429_1}
    &v(D,t) = 0, \\
    \label{430}
    &\dot{U}(t) = \phi(\vartheta,v,U,t),
    \end{align}   
with 
\begin{equation} \label{433}
\phi\left(\vartheta,v,U,t\right)\!=\!\! -cU(t)\!+\!c K\left[G(t)\!+\!\hat{H}(t)\!\!\int_{0}^{D}\!\!\!\!(D-r)(v(r,t)\!+\!U(t))dr\right].
\end{equation}
Next, we express the PDE system in \eqref{429}-\eqref{430} as an evolutionary equation in the Banach space \( \mathcal{X} \) \citep{27}, yielding $\dot{V}(t) = \mathcal{A}V(t) - \tilde{\phi}(\vartheta,V,U,t), t > 0$. Where \( V(t) \in \mathcal{X} \), and \( \mathcal{A} \) is the second-order derivative operator with Dirichlet and Neumann boundary conditions in \( \mathcal{X} \). To express \( v(0,t) \) in the ODE \eqref{428} in terms of \( V(t) \), we introduce the boundary operator \( \mathcal{B} \) and differential operator \( \mathcal{A} \)

\begin{align} \label{435}
\mathcal{A}\varphi := \frac{\partial ^2 \varphi}{\partial x^2}, \qquad \mathcal{B}V(t) := v(0,t),
\end{align}
and the domain
\begin{equation} \label{436}
D(\mathcal{A}) = \left\lbrace \varphi \in \mathcal{X} : \varphi, \frac{d}{dx} \varphi \in \mathcal{X} \text{ are a.c.}, \frac{d^2}{dx^2} \varphi \in \mathcal{X}, \frac{d}{dx} \varphi(0) = 0, \varphi(D) = 0 \right\rbrace.
\end{equation}    
Furthermore, we define the linear operators $\alpha^\top : \mathcal{X}\rightarrow \mathbb{R}$ and $\beta: \mathbb{R} \rightarrow \mathcal{X}$  as
\begin{align} \label{438}
\alpha^\top V(t) & := \int_{0}^{D}(D-r)v(r,t)\,dr, \qquad 
\beta \zeta := [\beta_1, \beta_2, \ldots]^\top \zeta, \qquad \zeta \in \mathbb{R},
\end{align} 
where \( \beta_k = -\sqrt{\frac{2}{D}} \frac{2D}{\pi (2k-1)} (-1)^k \),\( \psi_k(x) = \sqrt{\frac{2}{D}} \cos\left( \frac{\pi}{2} (2k-1) \frac{x}{D} \right) \), and \( k=1,2,\ldots\) are the eigenfunctions of \( \mathcal{A} \). Finally, the error dynamics with the infinite-dimensional state vector \( z(t) = [\vartheta(t) \, V(t) \, U(t)]^\top \) can be rewritten as \eqref{427} with
\begin{align} \label{439}
&\Gamma = \begin{bmatrix}
0 & \mathcal{B} & 1 \\
0 & \mathcal{A} & c\beta \\
0 & 0 & -c 
\end{bmatrix},\ J(\omega t,z)=\begin{bmatrix}
0 \\
-c\beta K\left[G(t)+\hat{H}(t)g(z)\right] \\
c K\left[G(t)+\hat{H}(t)g(z)\right]
\end{bmatrix},
\end{align} 
where \( g(z) = \frac{1}{2} D^2 U(t) + \alpha^\top V(t) \). Since \( \mathcal{A} \) generates an analytic semigroup \citep{28} and \( \mathcal{B} \) is \( \mathcal{A} \)-bounded (\( \|\mathcal{B} V(t)\| \leq 4D \sqrt{D} \|\mathcal{A} V(t)\| \)), the matrix \( \Gamma \) generates an analytic semigroup by the operator matrix theorem \citep{29}. Also, \( J(\omega t, z) \) in \eqref{439} is Fréchet differentiable in \( z \), continuous, and almost periodic in \( t \), uniformly with respect to \( z \). Thus, the error dynamics \eqref{42}-\eqref{45} have an exponentially stable periodic solution \( z^\Pi(t) \) satisfying \eqref{401}.

\subsubsection*{\textbf{Step 6:}\ \ \textit{Convergence to a Neighborhood of the Extremum}}
In this step, we prove the convergence statements \eqref{402}-\eqref{404}. leads to
\begin{align} \label{440}
\tilde{\theta}^2(t)\le 3\vartheta(t)^2+(4D+1)\|\bar{\alpha}_x\|^2.
\end{align}
Taking the time derivative of \( \bar{\Upsilon}(t) = \frac{\mu}{16D} \|\bar{\alpha}_x(t)\|^2 + \frac{1}{2} \vartheta^2(t) + \frac{1}{2} \|u(t)\|^2 + \frac{1}{2} \|u_x(t)\|^2 + \frac{1}{2} U(t)^2 \), and using the exponential stability of the original system, we get \( \dot{\bar{\Upsilon}} \le -\mu \min\left\{ \frac{1}{4}, \frac{1}{64D^3} \right\} \bar{\Upsilon} \), where \( \mu > 0 \). Thus, there exists some \( M > 0 \), such that
\begin{align} \label{441}
	\|\bar{\alpha}_x(t)\|^2 \le M e^{-k t}\ \ \text{with} \ \ k = \mu \min \left\lbrace \frac{1}{4},\frac{1}{64 D^3} \right\rbrace.
\end{align}
With \eqref{440}, \eqref{441} and adding the periodic solution $\vartheta^\Pi(t)$, it follows
\begin{align} \label{442}
\limsup\limits_{t\rightarrow \infty} |\tilde{\theta}(t)|^2 = \limsup\limits_{t\rightarrow\infty}\left\lbrace3|\vartheta(t)+\vartheta^\Pi(t)-\vartheta^\Pi(t)|^2\right\rbrace.
\end{align}
Using Young's inequality to \eqref{442}, we get 
\begin{align} \label{443}
|\vartheta(t) + \vartheta^\Pi(t) - \vartheta^\Pi(t)|^2 \leq \sqrt{2} \left( |\vartheta(t) - \vartheta^\Pi(t)|^2 + |\vartheta^\Pi(t)|^2 \right), 
\end{align}
which implies that \( \vartheta(t) - \vartheta^\Pi(t) \to 0 \) exponentially. Along with \eqref{401}, we get
\begin{align} \label{444}
\limsup\limits_{t \to \infty} |\tilde{\theta}(t)| = \mathcal{O}\left(1\right) \quad \text{with} \quad \limsup\limits_{t \to \infty} |\tilde{\theta}(t)|^2 = \limsup\limits_{t \to \infty} \left\lbrace 3\sqrt{2} |\vartheta^\Pi(t)|^2 \right\rbrace.
\end{align}
Since $\theta(t) - \Theta^* = \tilde{\theta}(t) + S(t)$ from \eqref{12} and Fig ~\ref{ESloop}, we finally get with \eqref{443} the next ultimate bound 
\begin{align} \label{445}
\limsup\limits_{t\rightarrow \infty} |\theta (t)-\Theta^*| = \mathcal{O}\left(|a|e^{D\sqrt{\omega/2}}+1/\omega\right).
\end{align}
The convergence of \( \Theta(t) \) to \( \Theta^* \) is established in the next step. From \eqref{13} and taking the absolute value, we obtain \( |\Theta(t) - \Theta^*| = \vartheta(t) + a \sin(\omega t) \). As with the convergence of \( \theta(t) \) to \( \Theta^* \), we express \( \vartheta(t) \) in terms of the periodic solution \( \vartheta^\Pi(t) \), apply Young's inequality, and observe that \( \vartheta(t) - \vartheta^\Pi(t) \) decays exponentially. Thus, combining with \eqref{401}, we get
\begin{align} \label{447}
\limsup\limits_{t\rightarrow \infty}\left|\Theta(t)-\Theta^*\right| = \mathcal{O}\left(|a|+1/\omega\right).
\end{align}
To show the convergence of the output $y(t)$ of the static map to the optimal value $y^*$, we replace $\Theta(t)-\Theta^*$ in \eqref{05} by \eqref{13} and take the absolute value
\begin{align} \label{448}
|y(t)-y^*| = \left|\frac{H}{2}\left[\vartheta(t)+a\sin(\omega (t))\right]^2\right| .
\end{align}
Expanding the quadratic term in \eqref{448} and applying Young's inequality, we get \( |y(t) - y^*| = \left| H \left[ \vartheta(t)^2 + a^2 \sin^2(\omega t) \right] \right| \). As before, adding the periodic solution \( \vartheta^\Pi(t) \), applying Young's inequality, and noting that \( \vartheta(t) - \vartheta^\Pi(t) \) decays exponentially (by the averaging theorem \citep{30}), we obtain \eqref{404} along with \eqref{401}. \hfill $\square$

\begin{remark}
Although the diffusion PDE has been adopted as the benchmark reference case in this study, it is crucial to emphasize that the proposed framework does not impose type-specific constraints on PDEs. As demonstrated by systematic numerical experiments in the next section, PINNs exhibit feasible motion planning capabilities in canonical PDE systems including  parabolic-type PDEs, hyperbolic-type PDEs, and PDEs with moving boundaries.  
\end{remark}

\section{Simulation results}

PDEs can be categorized into two main types: parabolic PDEs and hyperbolic PDEs. Parabolic PDEs are commonly used to model phenomena such as diffusion and heat conduction, while hyperbolic PDEs are typically applied to describe wave propagation and dynamic processes. Specifically, Distributed PDEs and RAD PDEs belong to the parabolic category, while Transport Hyperbolic PDEs and Wave PDEs fall under the hyperbolic category. PDEs with moving boundaries represent a special type of PDE, and many PDEs have the potential to exhibit moving boundaries in certain situations. In this paper, we will look at five distinct types of PDEs, focusing on their equation forms, boundary controllers, and trajectory generation, as outlined in Table \ref{tab:PDEs}. A unified PINN framework will be used to solve these PDEs, and the designed ES strategy will be applied to process the system composed of these five PDEs.

\begin{table}[h]
    \centering
    \caption{ES for distinct classes of PDE systems}\label{tab:PDEs}
    \label{tab:es_pde}
    \begin{tabular}{p{0.23\textwidth} p{0.74\textwidth}}
        \toprule
        \multicolumn{2}{l}{Distributed delay} \\
        & \textbf{PDE:} \(
        \partial_t \alpha(x,t) = \partial_x \alpha(x,t), \quad x \in [0,D], \quad y = Q \left( \int_0^D \Theta(t-\sigma) d\beta(\sigma) \right)
        \) \\
        & \textbf{BC (Dirichlet):} \(
        U(t) = \mathcal{T} \left\{ k \left[ G(t) + \hat{H}(t) \int_0^D (1-\beta(\sigma)) u(D-\sigma,t) \, d\sigma \right] \right\}
        \) \\
        & \textbf{Trajectory generation:} \(
        S(t) = \frac{a}{\gamma(\omega)} \int_0^D \sin(\omega (t + \xi)) \, d\beta(\xi).
        \) \\        
        \midrule
        \multicolumn{2}{l}{RAD equation} \\
        & \textbf{PDE:} \(
        \partial_t \alpha(x,t) = \epsilon \partial_{xx} \alpha(x,t) + b \partial_x \alpha(x,t) + \lambda \alpha(x,t), \quad x \in [0,1]\) \\
        & \textbf{BC (Dirichlet):} \(
        U(t) = \mathcal{T} \{ k e^{-\frac{b}{2\epsilon}} [ \gamma(1) G(t) + \hat{H}(t) \int_0^1 e^{\frac{b}{2\epsilon} \sigma} m(1-\sigma) \} \) \( \times u(\sigma,t) \, d\sigma ] \)
        \(\gamma(x) = \cosh \left( \sqrt{\frac{\xi}{\epsilon}} x \right) + \frac{b}{2\epsilon} \sqrt{\frac{\epsilon}{\xi}} \sinh \left( \sqrt{\frac{\xi}{\epsilon}} x \right)\), \(\xi := \frac{b^2}{4\epsilon} - \lambda \leq 0\), \(k > 0\), \(m(x-\sigma) = \frac{1}{\epsilon} \sqrt{\frac{\xi}{\epsilon}} \sinh \left( \sqrt{\frac{\xi}{\epsilon}} (x-\sigma) \right)\), \(\epsilon > 0\), \(b \geq 0\), \(\lambda \geq 0\). \\
        & \textbf{Trajectory generation:} \(
        S(t) = e^{\frac{b}{2\epsilon}} \sum_{k=0}^{\infty} \frac{a_{2k}(t)}{(2k)!} + \frac{b}{2\epsilon} \frac{a_{2k}(t)}{(2k+1)!}, \)
        \(a_{2k} := \frac{a}{\epsilon^k} \sin(\omega t) \sum_{n=0}^k \binom{k}{2n} \xi^{k-2n} \omega^{2n} + \frac{a}{\epsilon^k} \cos(\omega t) \sum_{n=0}^k \binom{k}{2n+1} \xi^{k-2n-1} \omega^{2n+1}\). \\
        \midrule
        \multicolumn{2}{l}{Transport PDE} \\
        & \textbf{PDE:} \(
        \partial_t u(x,t) + c \partial_x u(x,t) = D u(x,t - D)
        \), \( \quad D \geq 0 \), \( \quad x \in [0,D] \) \\
        & \textbf{BC (Neumann):} \(
        U(t) = \mathcal{T} \left\{ k \left[ G(t) + \hat{H}(t) \int_{t-D}^{t} U(\tau) \, d\tau \right] \right\}
        \) \\
        & \textbf{Trajectory generation:} \(
        S(t) = a \sin(\omega (t+D)).
        \) \\
        \midrule
        \multicolumn{2}{l}{Wave dynamics} \\
        & \textbf{PDE:} \(
        \partial_{tt} \alpha(x,t) = \partial_{xx} \alpha(x,t), \quad x \in [0,D]
        \) \\
        & \textbf{BC (Neumann):} \(
        U(t) = \mathcal{T} \{ c \left[ k \hat{H}(t) u(D,t) - \partial_t u(D,t) \right] + \rho(D) G(t) \), \(+ \hat{H}(t) \int_0^D \rho(\sigma) \partial_t u(\sigma,t) \, d\sigma \Big\},\)
        \(\rho(s) = k \begin{bmatrix} 0 & \Gamma_1 e^{As} \end{bmatrix} \begin{bmatrix} I \\ A \end{bmatrix}\), \(A = \begin{pmatrix} 0 & 0 \\ I & 0 \end{pmatrix}\). \\
        & \textbf{Trajectory generation:} \(
        S(t) = a \cos(\omega D) \sin(\omega t).
        \) \\
        \midrule
        \multicolumn{2}{l}{Stefan PDE:} \\
        & \textbf{PDE} \(
        \partial_t \alpha(x,t) = \partial_{xx} \alpha(x,t), \quad x \in (0, s(t))., \quad x \in [0,D]
        \) \\
        & \textbf{BC (Dirichlet):} \(
        U(t) = \mathcal{T} \left\{ K \left[ G(t) + \hat{H}(t) \int_0^{s(t)} u(x,t) \, dx \right] \right\}
        \) \\
        & \textbf{Trajectory generation:} \(
        S(t) = - \sum_{i=0}^{\infty} \frac{1}{(2i - 1)!} \frac{\partial^i}{\partial t^i} \left[ -a \sin(\omega t) \right]^{2i - 1}
        \) \\
        \bottomrule
    \end{tabular}
\end{table}

\subsection{Comparison of PINN-Based and Numerical Solution Results}

First, a unified PINN framework was utilized to solve the five partial differential equations listed in Table~\ref{tab:PDEs}. The framework, leveraging the powerful expressive capacity of neural networks, can effectively handle partial differential equations with complex boundary and initial conditions. The network consisted of 5 layers, each containing 100 neurons. The activation function used in each layer was Tanh, and the optimizer was selected as the Adam optimization algorithm, which enabled the network to efficiently converge to the optimal solution. Table~\ref{tab:network_parameters} provides a detailed overview of these key network parameters. Additionally, the boundary condition \( u(0,t) = \sin(5\omega t) \) was applied to the system to assess the framework's performance under various conditions.

\begin{figure}  
    \centering
    \subfloat[Numerical solution results of Distributed PDE]{
        \includegraphics[width=0.29\linewidth]{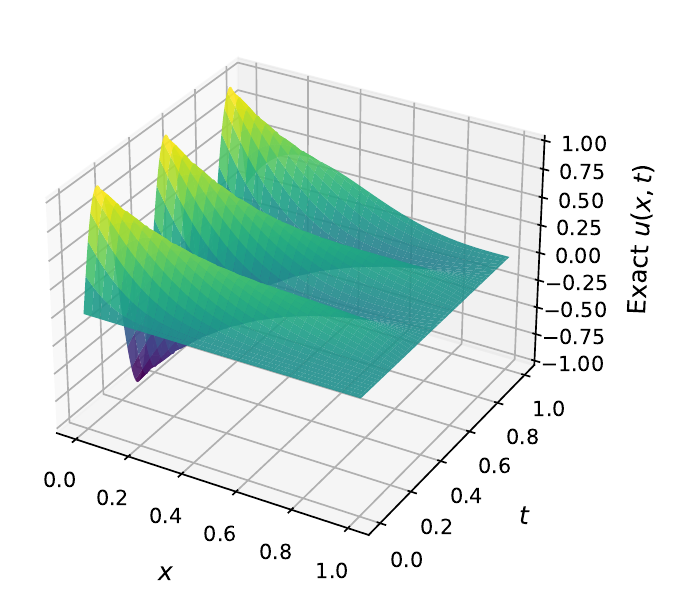}
    }
    \hspace{0.25cm}  
    \subfloat[PINN solution results of Distributed PDE]{
        \includegraphics[width=0.29\linewidth]{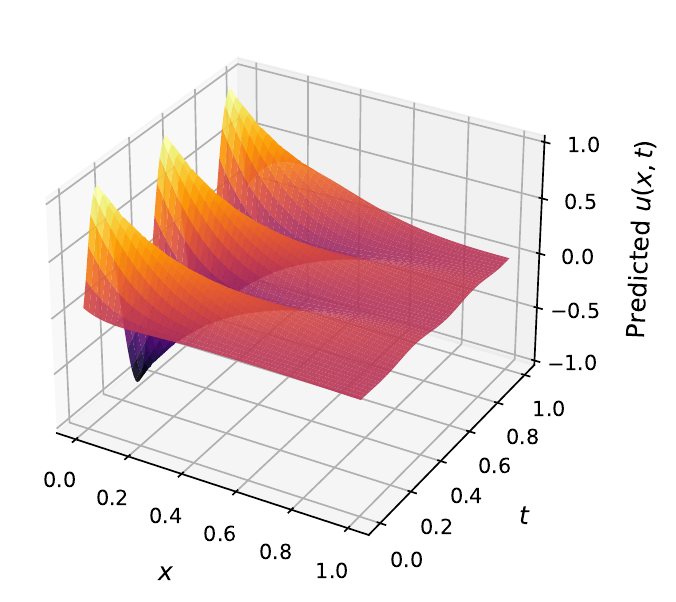}
    }
    \hspace{0.25cm}  
    \subfloat[Distributed PDE: Absolute Error]{
        \includegraphics[width=0.29\linewidth]{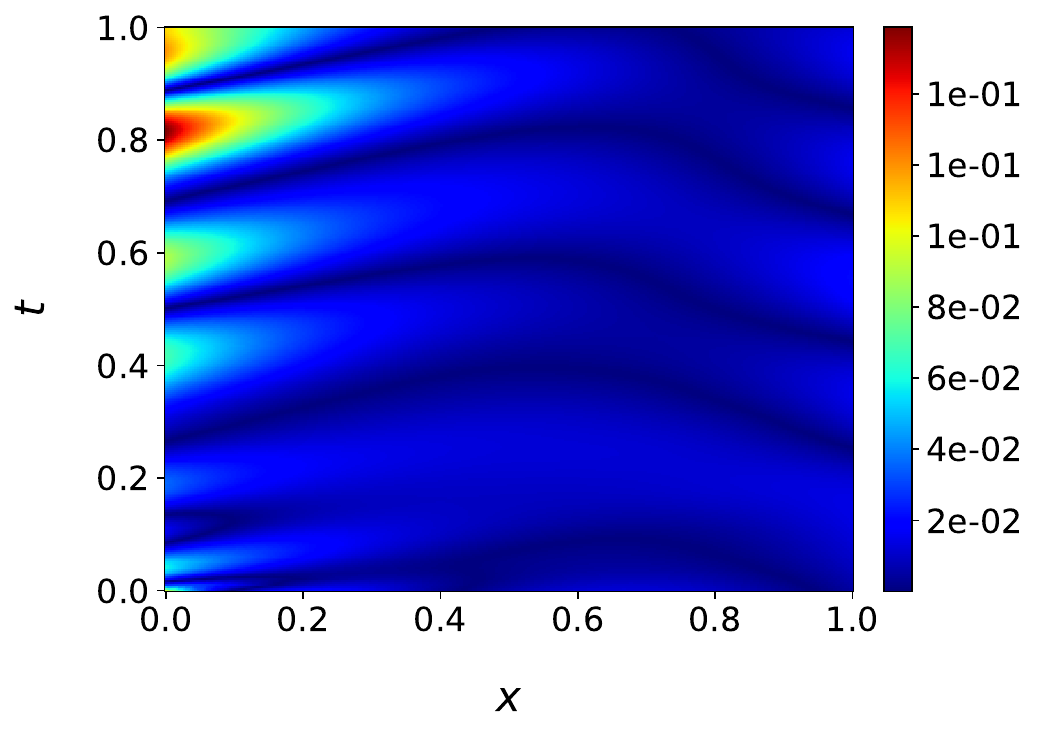}
    }
    \vspace{-0.5cm}  
    \subfloat[Numerical solution results of RAD PDE]{
        \includegraphics[width=0.29\linewidth]{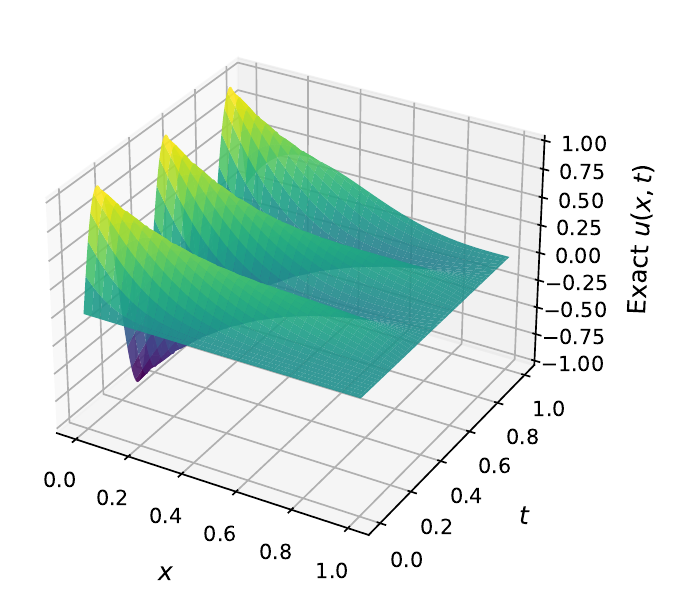}
    }
    \hspace{0.25cm}  
    \subfloat[PINN solution results of RAD PDE]{
        \includegraphics[width=0.29\linewidth]{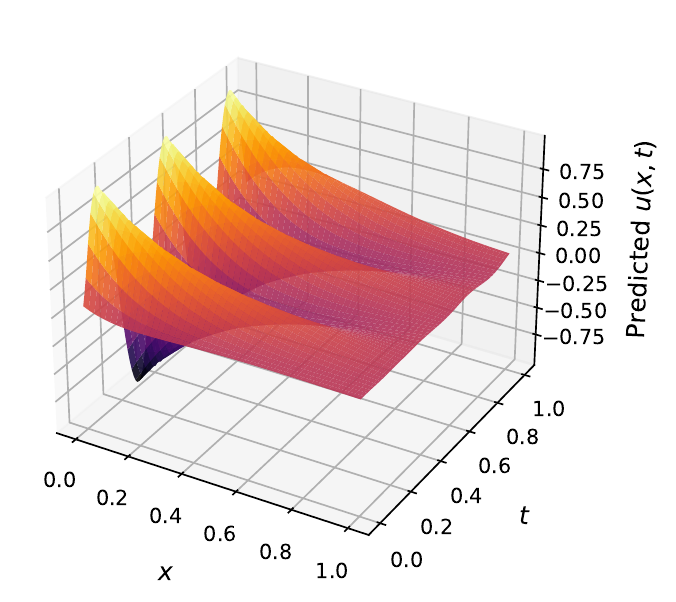} 
    }
    \hspace{0.25cm}  
    \subfloat[RAD PDE: Absolute Error]{
        \includegraphics[width=0.29\linewidth]{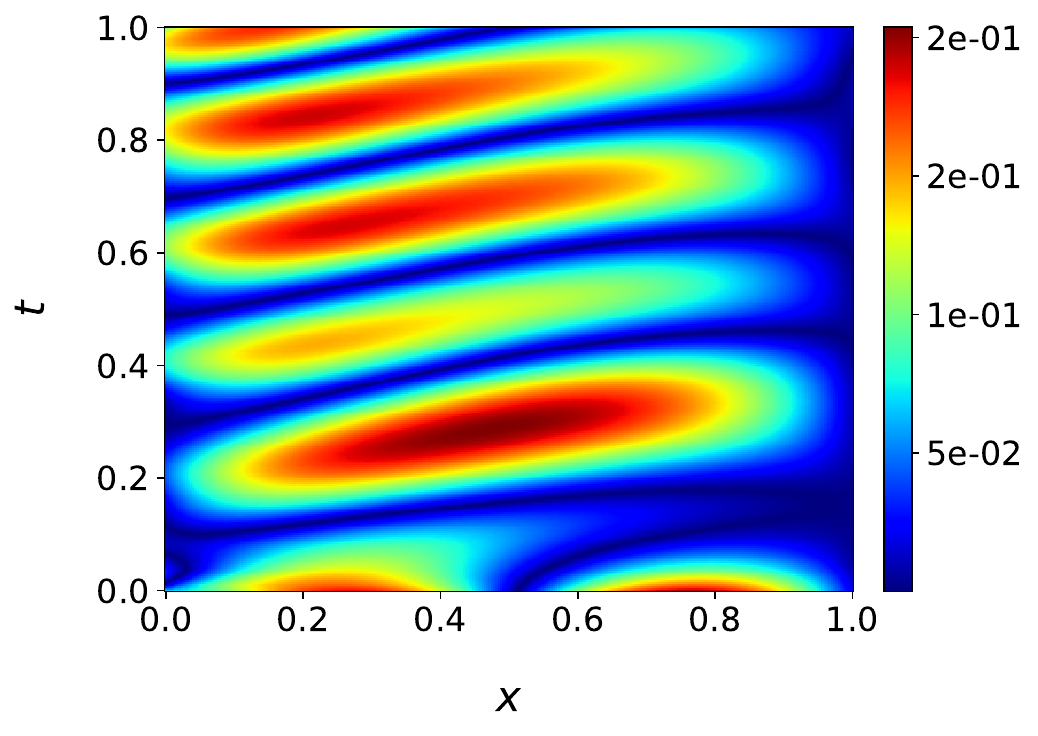} 
    }
    \vspace{-0.5cm}  
    \subfloat[Numerical solution results of Wave PDE]{
        \includegraphics[width=0.29\linewidth]{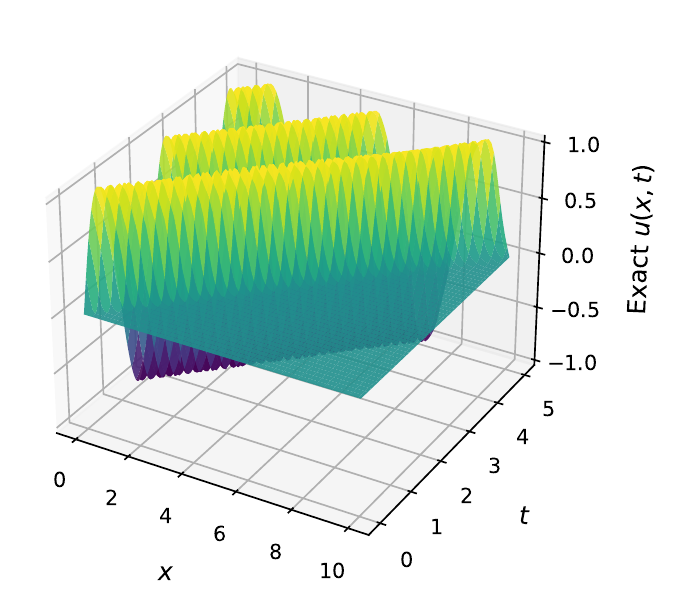}
    }
    \hspace{0.25cm}  
    \subfloat[PINN solution results of Wave PDE]{
        \includegraphics[width=0.29\linewidth]{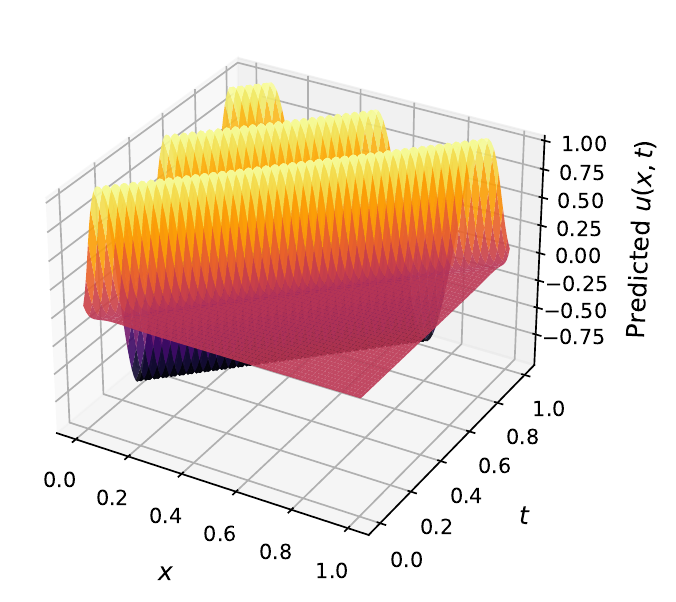} 
    }
    \hspace{0.25cm}  
    \subfloat[Wave PDE: Absolute Error]{
        \includegraphics[width=0.29\linewidth]{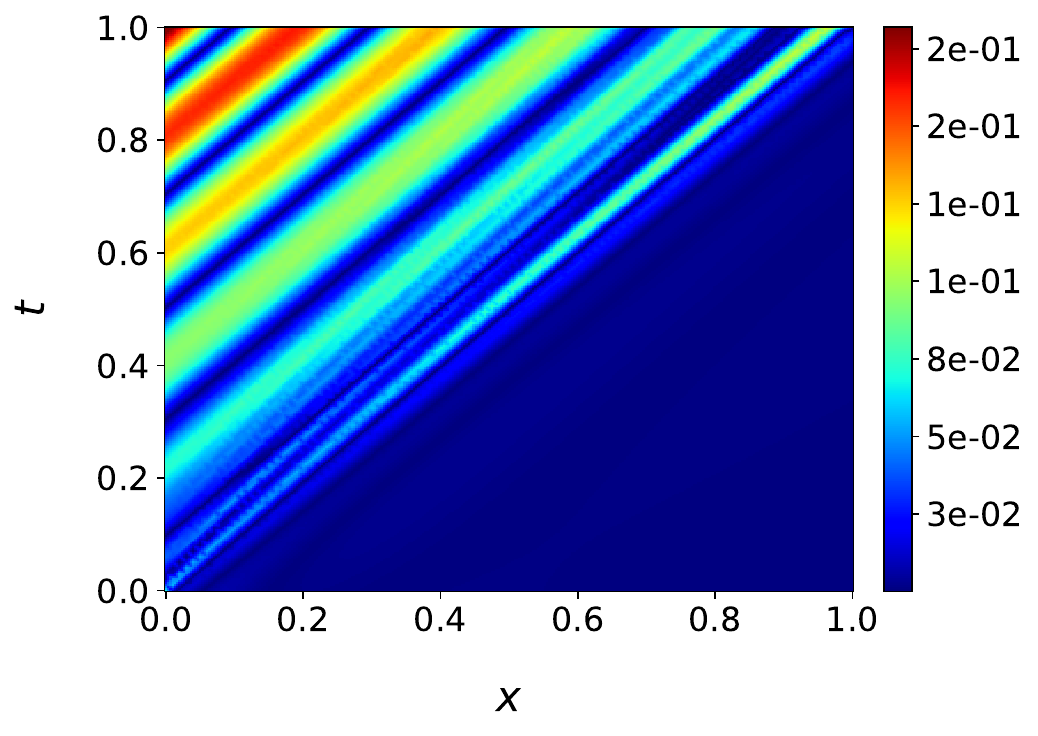} 
    }
    \vspace{-0.5cm}  
    \subfloat[Numerical solution results of Transport PDE]{
        \includegraphics[width=0.29\linewidth]{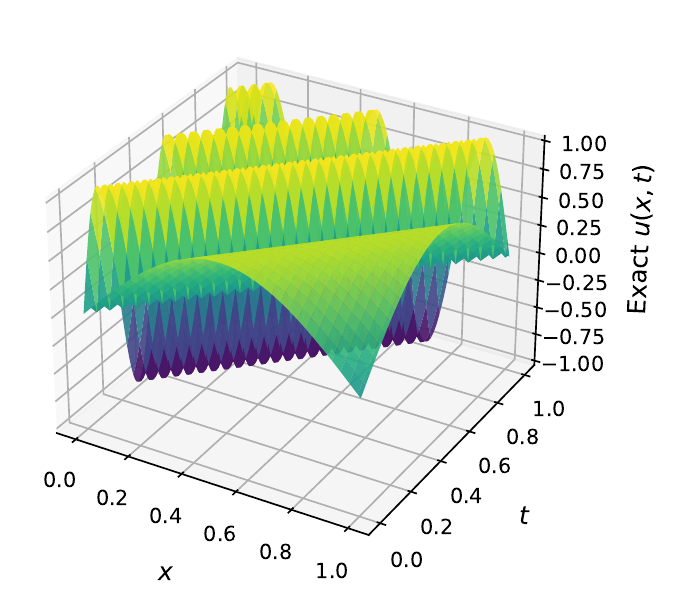}
    }
    \hspace{0.25cm}  
    \subfloat[PINN solution results of Transport PDE]{
        \includegraphics[width=0.29\linewidth]{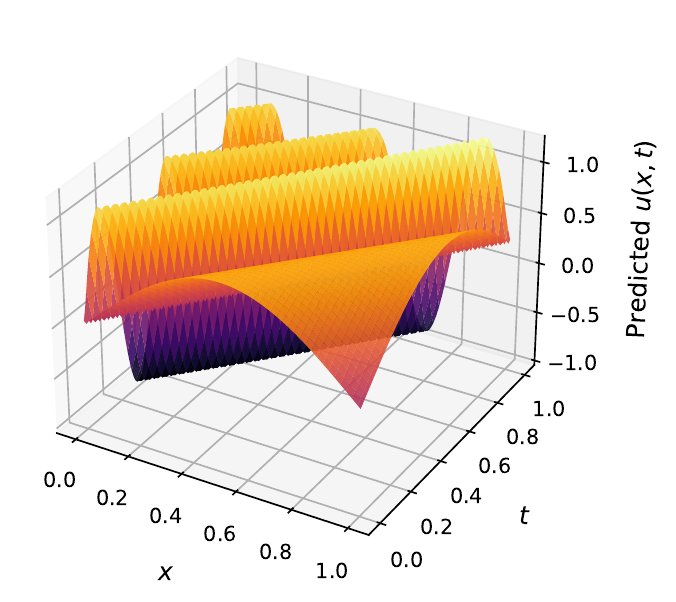} 
    }
    \hspace{0.25cm}  
    \subfloat[Transport PDE: Absolute Error]{
        \includegraphics[width=0.29\linewidth]{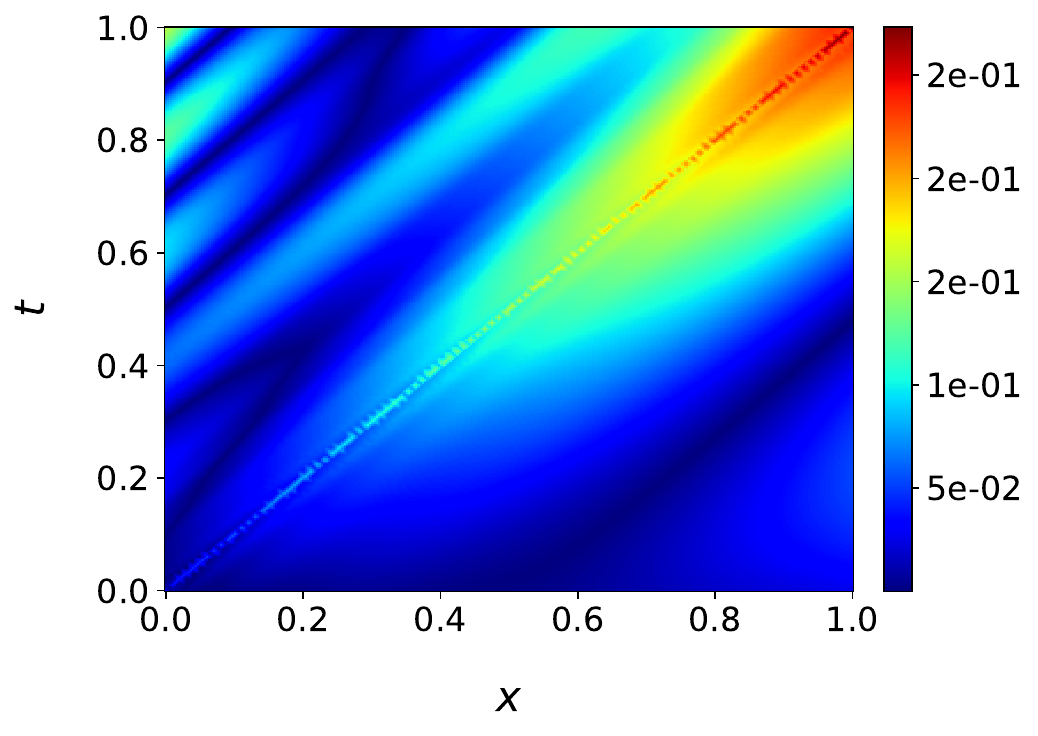} 
    }
    \caption{Comparison of Numerical and PINN Solutions for PDEs}
    \label{fig:PDEs}
\end{figure}

\begin{figure*}[htbp]  
    \centering
    \subfloat[Numerical solution results of Stefan PDE]{
        \includegraphics[width=0.29\linewidth]{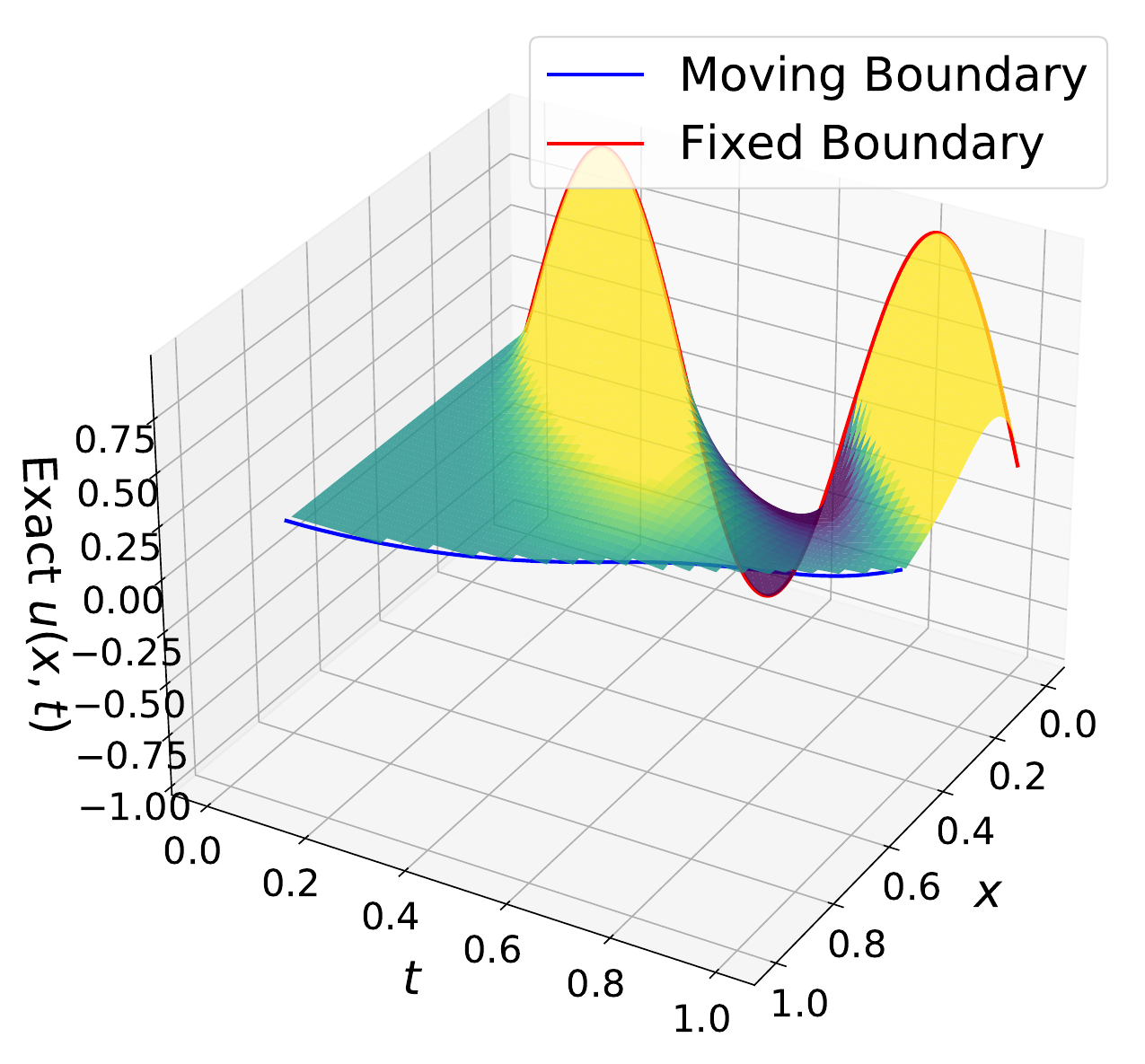}
    }
    \hspace{0.25cm}  
    \subfloat[PINN solution results of Stefan PDE]{
        \includegraphics[width=0.29\linewidth]{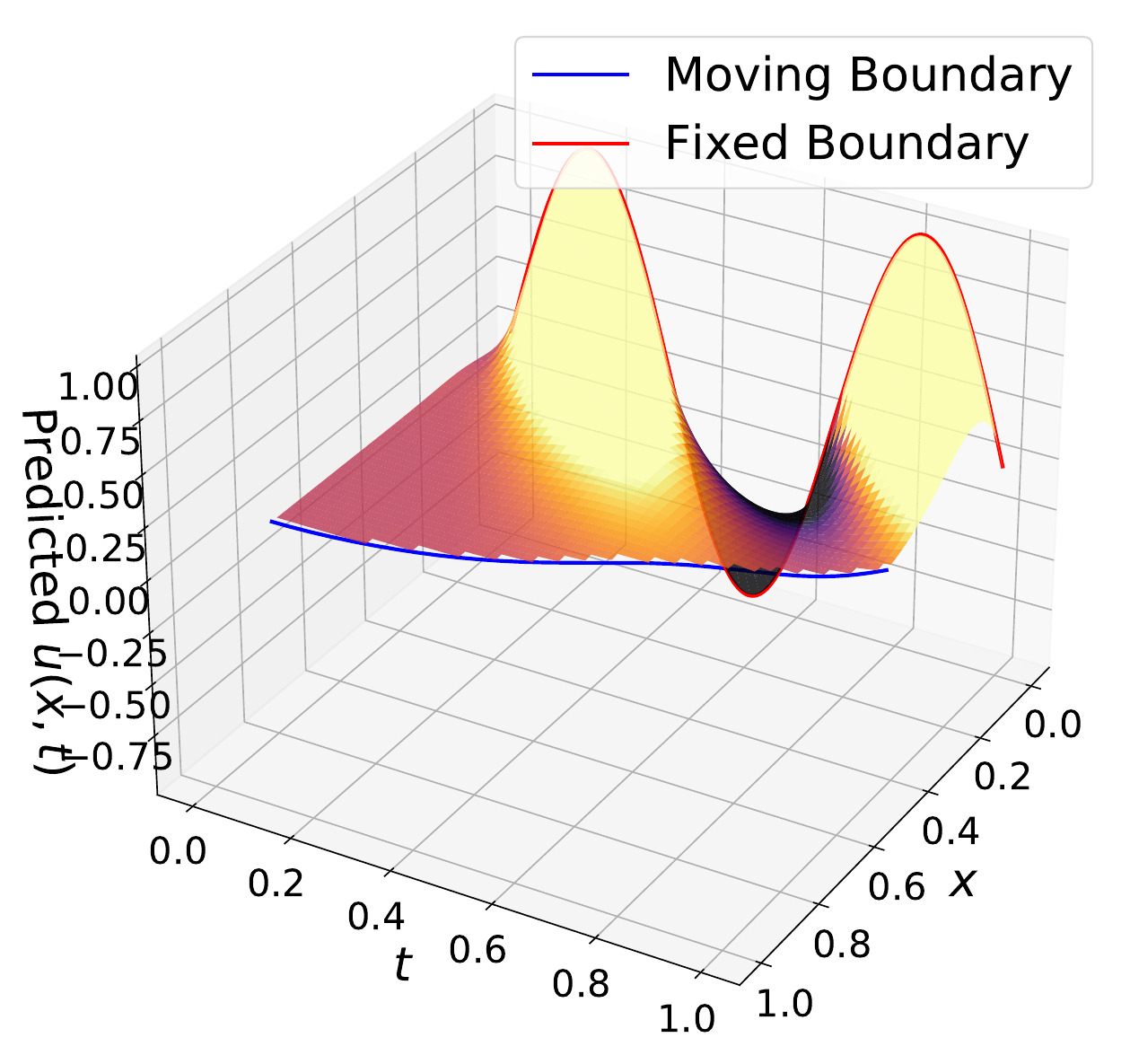}
    }
    \hspace{0.25cm}  
    \subfloat[Stefan PDE: Absolute Error]{
        \includegraphics[width=0.29\linewidth]{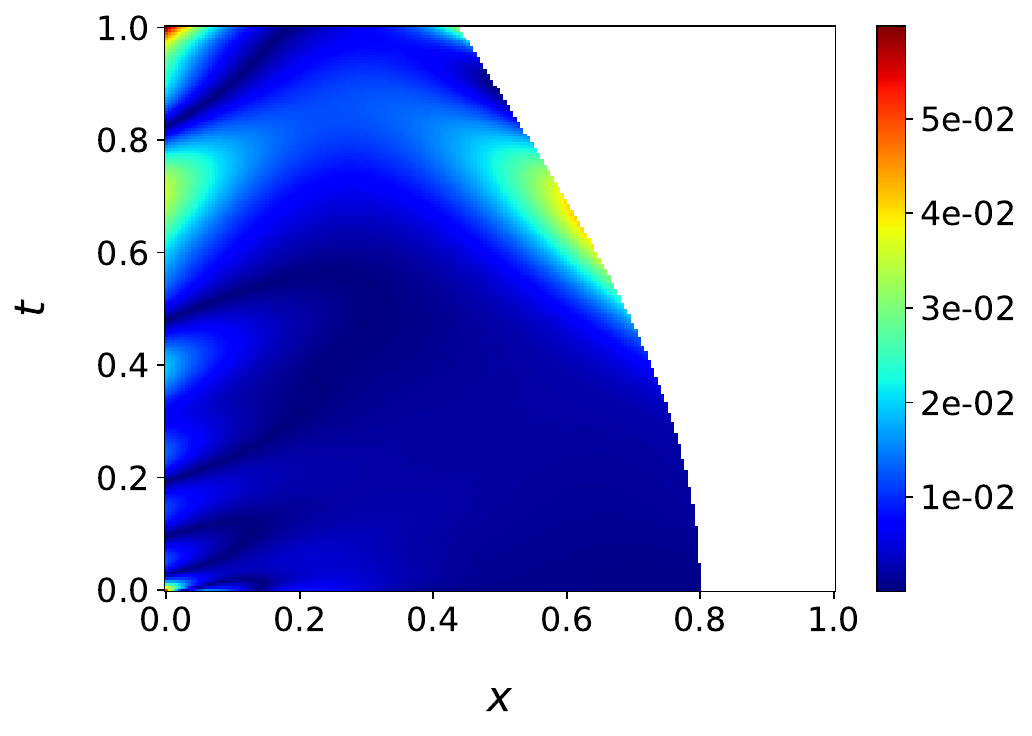}
    }
    \caption{Comparison of Numerical and PINN Solutions for Stefan PDE}
    \label{fig:Stefan_PDE}
\end{figure*}

\begin{table}
    \caption{Network Parameters\label{tab:network_parameters}}%
    \centering
    \begin{tabular*}{\columnwidth}{@{\extracolsep{\fill}} l l l l @{\extracolsep{\fill}}}
    \toprule
    \textbf{Parameter}                & \textbf{Value}     & \textbf{Parameter}          & \textbf{Value} \\
    \midrule
    Hidden Layers                     & 5                  & Neurons per Layer           & 100            \\
    Activation Function               & Tanh               & Optimizer                   & Adam           \\
    Initial Learning Rate             & 0.001              & Training Iterations          & 40,000         \\
    Training Data Points              & 10,000             & Programming Framework        & Tensorflow     \\
    \botrule
    \end{tabular*}
    \end{table}

In Fig.~\ref{fig:PDEs}, we present twelve images comparing the numerical solutions and PINN solutions for four types of PDEs: two diffusion-type and two parabolic-type. These include Distributed PDEs, RAD PDEs, Wave PDEs, and Transport PDEs. Each row contains three subplots: the numerical solution results, the PINN solution results, and the absolute error between the two methods. The results show that the PINN framework effectively captures the physical characteristics of all four PDE types. For diffusion-type PDEs, PINN accurately models the diffusion process, even under complex boundary conditions. For parabolic PDEs, PINN captures the time-evolution of the system and effectively handles both transient and steady-state behaviors. In Fig.~\ref{fig:Stefan_PDE}, we compare numerical methods and PINN methods for solving the Stefan PDE. The results demonstrate that PINNs are capable of handling time-varying boundary conditions effectively. The solutions obtained by PINNs closely match the numerical results, highlighting the framework's ability to capture dynamic system behaviors. This showcases the versatility and robustness of PINNs in solving the Stefan PDE, even in the presence of complex boundary interactions.

In Fig~\ref{fig:PDEs} (c) (f) (i) (l) and Fig~\ref{fig:Stefan_PDE} (c), we present the absolute errors between the PINN and numerical solutions for five PDEs. The figure highlights the distribution of these errors, revealing distinct patterns, though the error magnitudes vary significantly across equations. This variation is likely due to factors such as nonlinearity, boundary conditions, and computational precision. Overall, the error magnitudes are on the order of \(10^{-3}\), demonstrating the effectiveness and accuracy of the unified PINN framework in solving these PDEs.

\subsection{Simulation Results of Probing Signals}\label{sec:probing}
 
\begin{figure}[htbp] 
    \centering
    \begin{minipage}{0.45\textwidth}
        \centering
        \includegraphics[width=\textwidth]{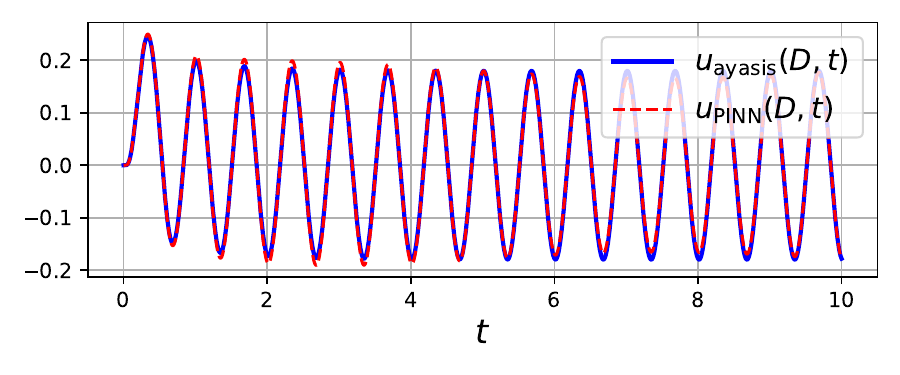}
        \caption{Distributed PDE perturbation signal evolution}  
        \label{p1}
    \end{minipage}
    \hfill  
    \begin{minipage}{0.45\textwidth}
        \centering
        \includegraphics[width=\textwidth]{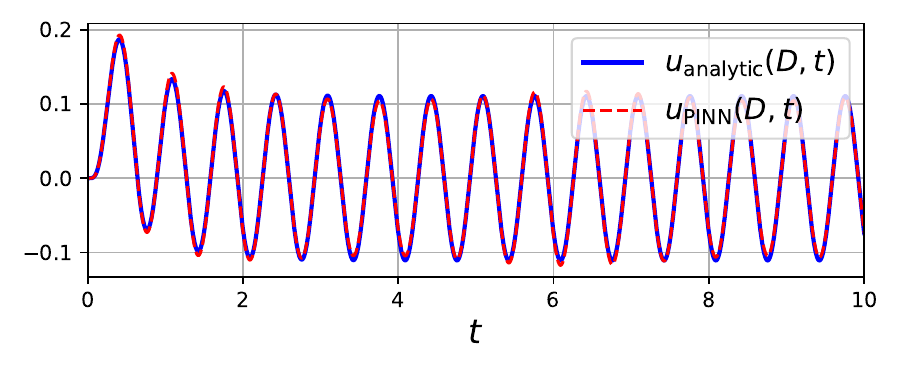}
        \caption{RAD PDE perturbation signal evolution}  
        \label{p2}
    \end{minipage}  
    \vspace{1em}  
    \begin{minipage}{0.45\textwidth}
        \centering
        \includegraphics[width=\textwidth]{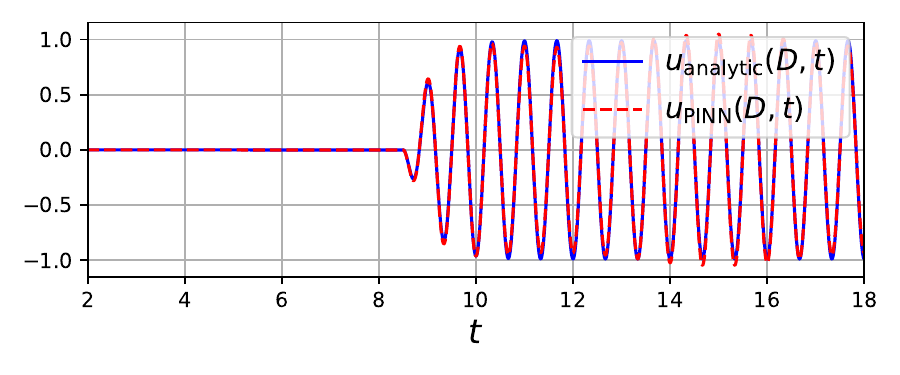}
        \caption{Wave PDE perturbation signal evolution}  
        \label{p3}
    \end{minipage}
    \vspace{-0.4cm}  
    \hfill  
    \begin{minipage}{0.45\textwidth}
        \centering
        \includegraphics[width=\textwidth]{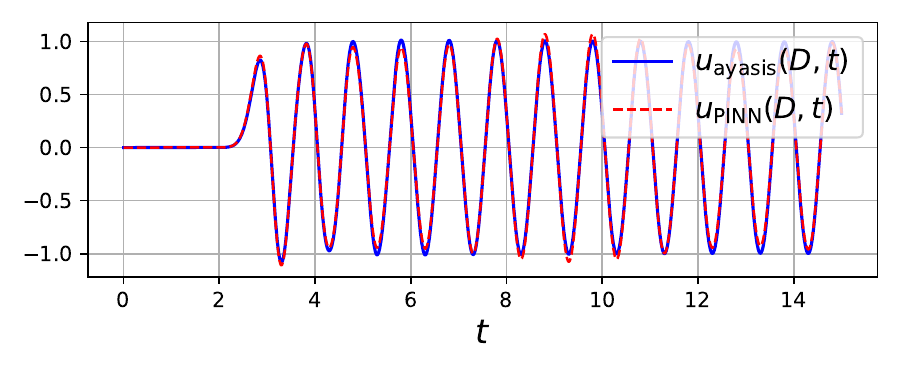}
        \caption{Transport PDE perturbation signal evolution}  
        \label{p4}
    \end{minipage}
    \vspace{-0.4cm}  
\end{figure}

We utilize the designed PINN framework to solve the perturbation signals of the five equations in Table \ref{tab:PDEs}, with identical left boundary conditions applied. Figures~\ref{p1} to \ref{p4} compare the temporal evolution of the PINN-predicted signals against analytical benchmarks.
The results show excellent agreement between PINN predictions and analytical solutions. For the distributed PDE, the dual exponential modulation is replicated with minimal error. In the RAD PDE, the network resolves both phase shift and amplitude decay driven by advection-diffusion coupling. The wave PDE solution demonstrates precise amplitude modulation at critical domain lengths, avoiding spurious oscillations typical of grid-based methods. For the transport PDE, the PINN captures the delayed instability onset without numerical dissipation. These observations highlight the outstanding capability of PINNs in handling diverse PDE types while preserving physical constraints.
\begin{figure}[htbp] 
    \centering
    \begin{minipage}{0.45\textwidth}
        \centering
        \includegraphics[width=\textwidth]{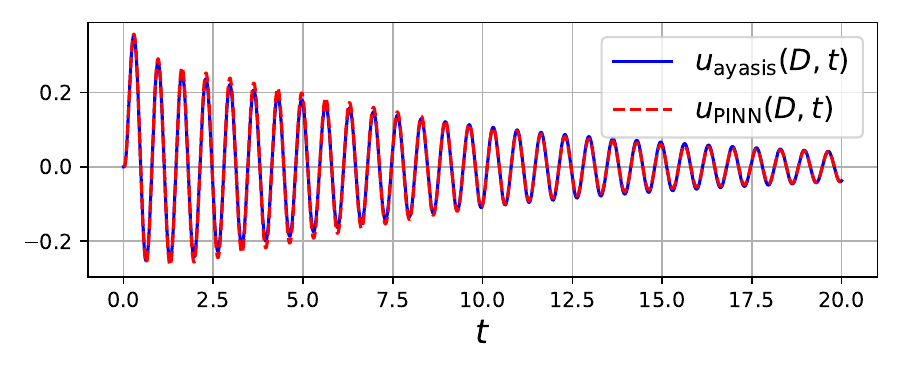}
        \caption{Perturbation signal evolution via Stefan PDE}
        \label{p5}
    \end{minipage}
    \vspace{-0.6cm}  
\end{figure}

For Stefan-type PDEs with a moving boundary \( s(t) = k_t e^{-k_b t} \), the perturbation signal \( S_S(t) \) must account for the time-varying spatial domain \( x \in [0, s(t)] \), where \( s(t) \) governs the domain length. The error order can be analyzed by the remainder term of the truncated series, where the size of the remainder is determined by the behavior of higher-order derivatives and truncation order, typically expressed as \( O(k^{n+1}) \), reflecting the rapid decay or growth of the error as the truncation order increases.

Compared to traditional numerical methods, PINNs approximate the solution through nonlinear superposition, expressed as \(\beta(x, t) \approx \sum_{i=1}^N c_i \sigma(w_i x + b_i t + a_i),\) where \( c_i \in \mathbb{R} \) are trainable weight coefficients and \( \sigma(\cdot) \) is a smooth activation function. The accuracy of PINNs stems from adaptive learning rate optimization and deep network architecture. A deep network with at least three hidden layers can simultaneously capture both exponential decay and oscillatory perturbations, as shown in Fig \ref{p5}. The numerical experiments in this paper show that for Stefan PDEs, the minimum truncation error in traditional numerical methods occurs at a truncation order of \( 10^{-2} \), with an error order of \( O(k^{n+1}) \). The PINN framework used to solve Stefan PDEs reduces the error to \( 10^{-4} \), outperforming the numerical methods by two orders of magnitude. PINNs demonstrate more precise results in solving Stefan PDEs and confirm their superiority in solving Stefan-type PDEs. This superiority arises because PINNs avoid the rigidity of polynomial expansions, instead resolving the implicit coupling between \( s(t) \) and diffusion dynamics through global optimization.

The above experiments demonstrate that PINNs can directly generate accurate perturbation signals, eliminating the need for trajectory generation in traditional numerical methods. Compared to the method of generating perturbation signal trajectories through numerical methods for five PDEs, as shown in Table \ref{tab:PDEs}, we find that PINNs not only align closely with the analytical benchmark in terms of time evolution, but also avoid common numerical errors and pseudo-oscillations in traditional methods through global optimization.

\begin{figure}[htbp]  
    \vspace{-0.6cm}
    \centering
    \subfloat[$\eta=8^{-3}$]{
        \includegraphics[width=0.24\linewidth]{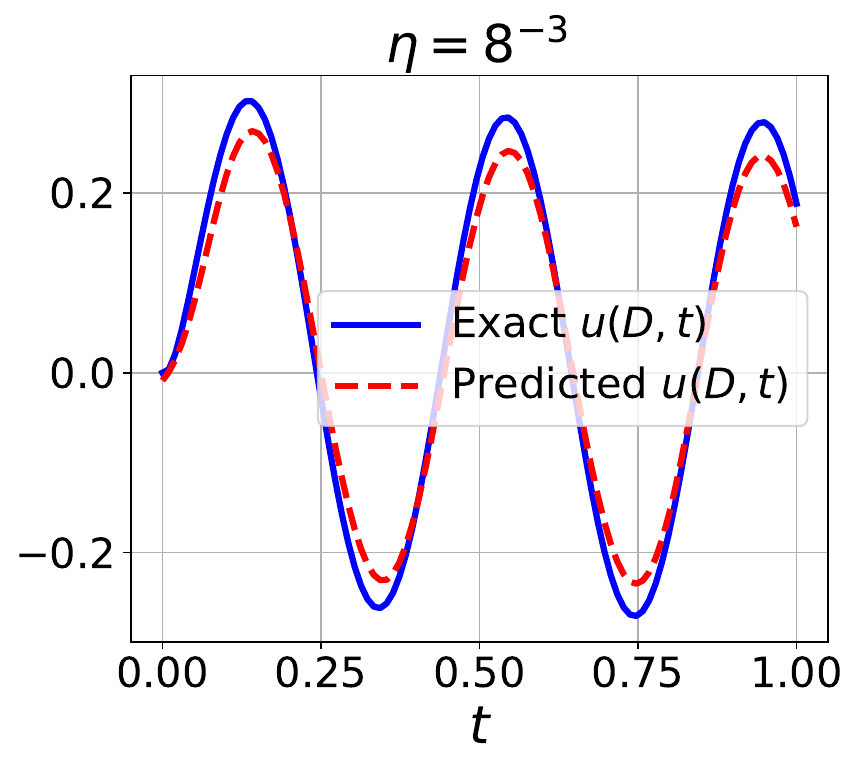}
    }
    \subfloat[$\eta=10^{-3}$]{
        \includegraphics[width=0.24\linewidth]{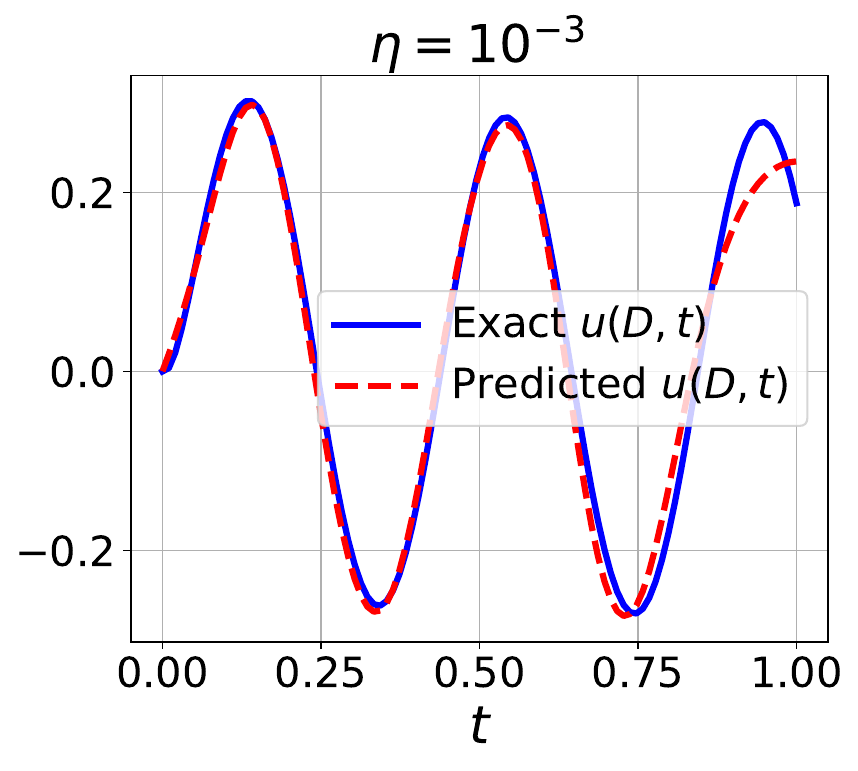}
    }
    \subfloat[$\eta=10^{-4}$]{
        \includegraphics[width=0.24\linewidth]{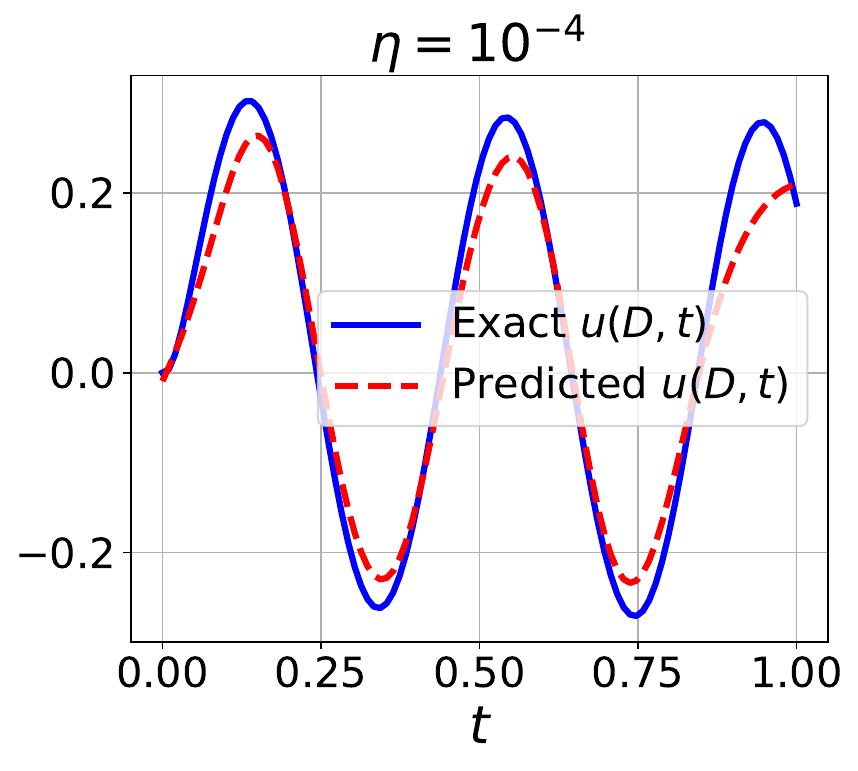}
    }
    \subfloat[$MSE vs. \eta$]{
        \includegraphics[width=0.24\linewidth]{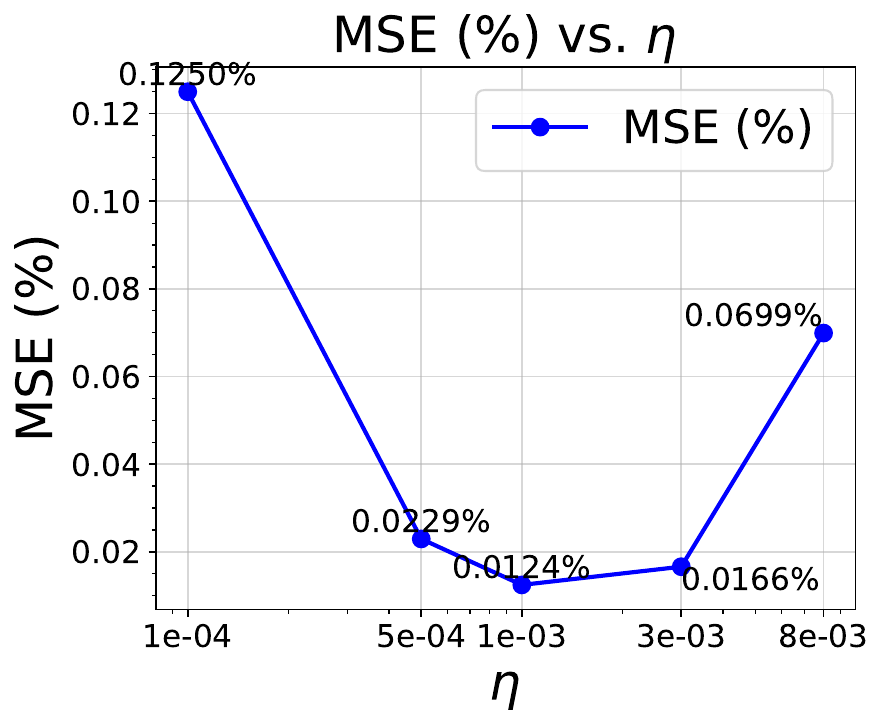}
    }
    \caption{PINN Approximation of MSE vs. Variation of $\eta$}
    \label{a_MSE}
    \vspace{-0.6cm}
\end{figure}

\begin{figure}[htbp]  
    \vspace{-0.6cm}
    \centering
    \subfloat[$batch\,size=32$]{
        \includegraphics[width=0.24\linewidth]{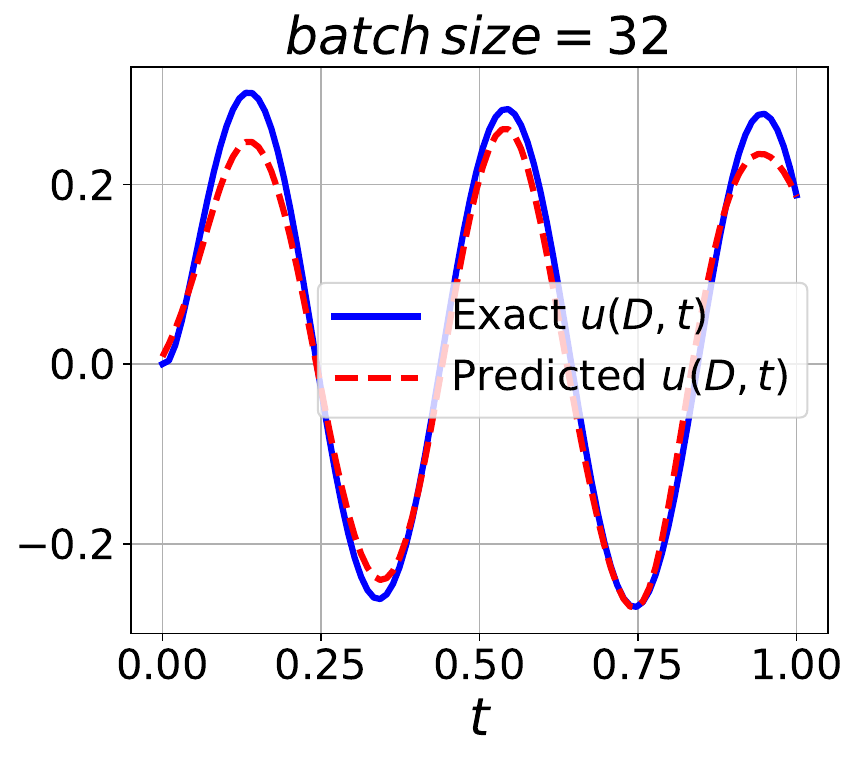}
    }
    \subfloat[$batch\,size=128$]{
        \includegraphics[width=0.24\linewidth]{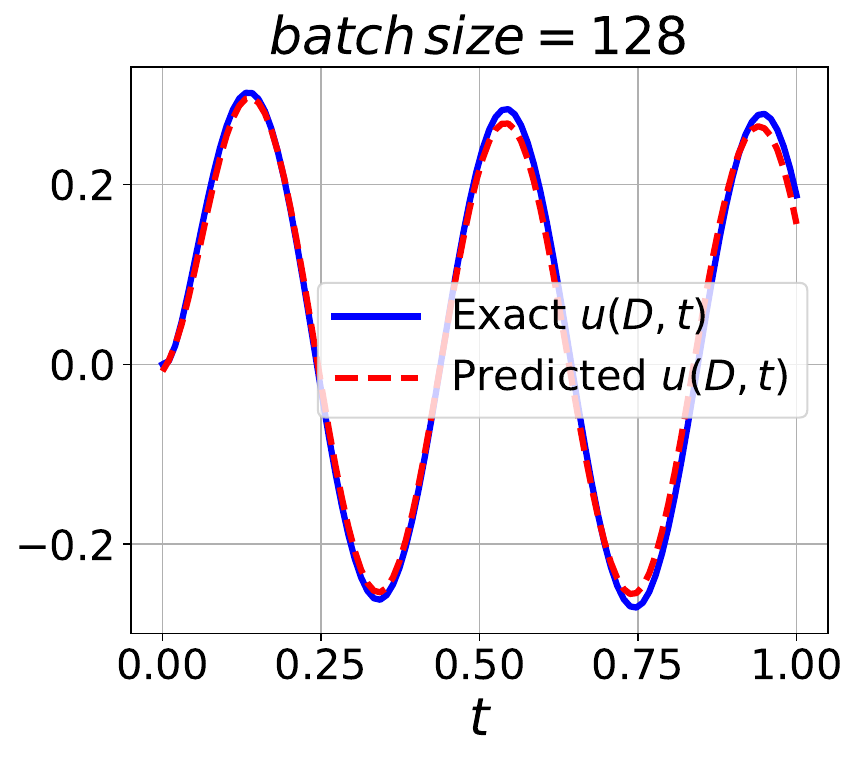}
    }
    \subfloat[$batch\,size=512$]{
        \includegraphics[width=0.24\linewidth]{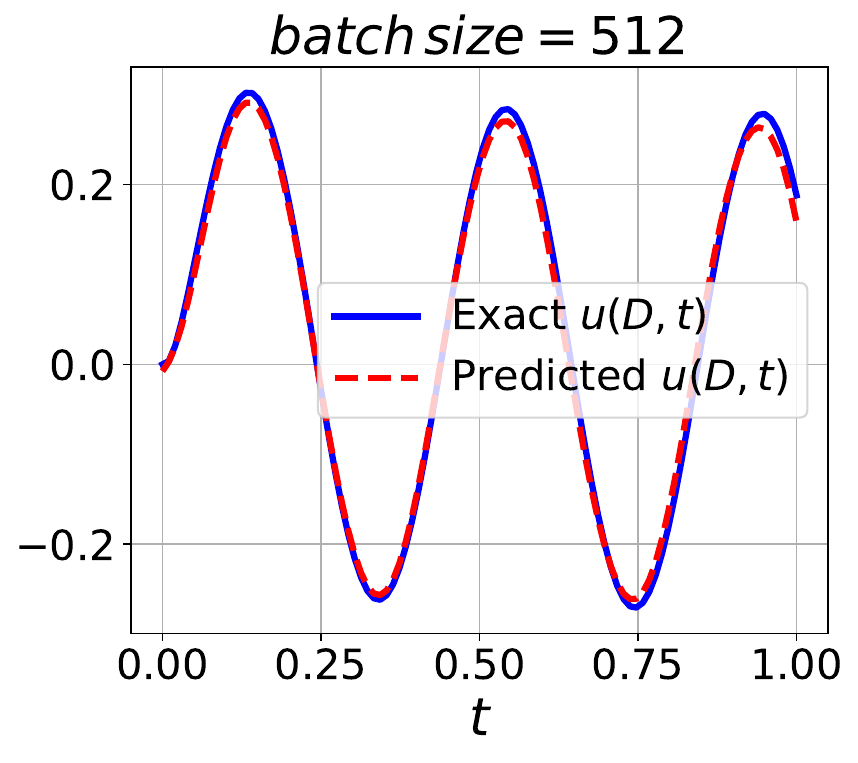}
    }
    \subfloat[$MSE vs. batch\,size$]{
        \includegraphics[width=0.24\linewidth]{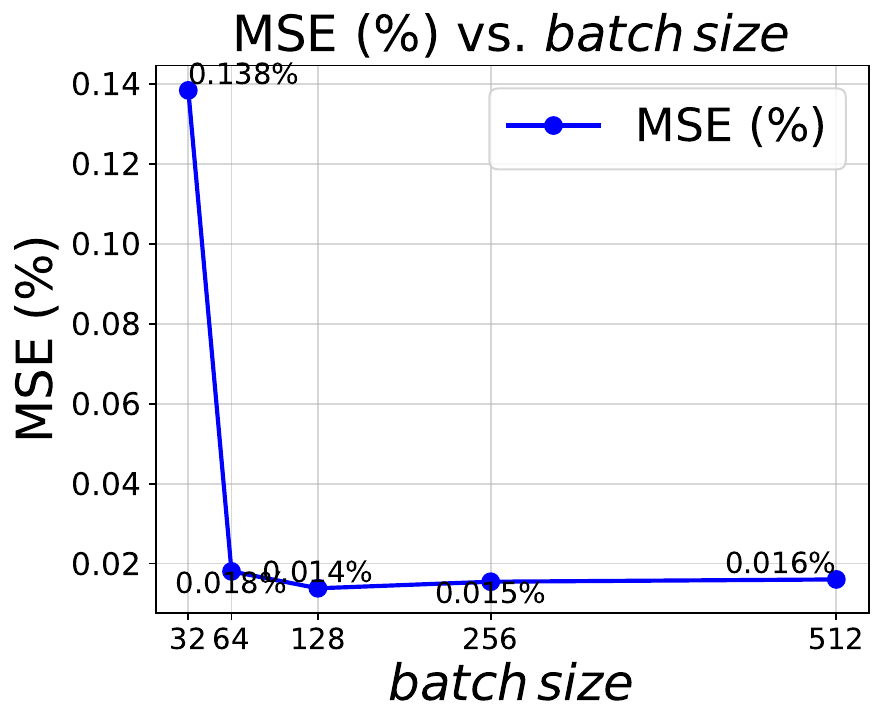}
    }
    \caption{PINN Approximation of MSE vs. Variation of $batch\,size$}
    \label{b_MSE}
    \vspace{-0.6cm}
\end{figure}

When generating \( S(t) \) through \( \text{asin}(wt) \), different PINN parameter settings have varying impacts on the neural network's approximation results. The most significant parameters are the learning rate and batch size. Using the diffusion PDE as an example, we analyze the effects of these two parameters on the perturbation signal \( S(t) \). The approximation performance is quantified by calculating the mean squared error (MSE). Figure 10 illustrates the impact of different learning rates on the PINN approximation of \( S(t) \). From Fig~\ref{a_MSE}(d), it is evident that the best approximation occurs when the learning rate is approximately \( \eta = 10^{-3} \). Figure 11 shows the effect of different batch sizes on the PINN approximation of \( S(t) \). From Fig~\ref{b_MSE}(d), it can be observed that when the batch size is less than 128, the MSE decreases as the batch size increases. However, when the batch size exceeds 128, the approximation performance remains relatively unchanged. Moreover, larger batch sizes result in longer training times. Therefore, a batch size of around 128 provides a good balance for achieving optimal approximation performance.

\subsection{Simulation Results of Extremum Seeking Control}

In this study, we employ an extremum seeking system for five types of PDEs. Specifically, the system is configured as follows: First, consider actuation dynamics governed by diffusion PDEs, where the Hessian matrix is \( H = -2 \), the optimizer is \( \Theta^* = 2 \), the optimal value is \( y^* = 5 \), and the domain length for the actuator dynamics is \( D = 1 \). The parameters of the dithering signal and the designed controller are set as \( \omega = 10 \), \( a = 0.2 \), \( c = 10 \), and \( K = 0.2 \). Second, the perturbation signals are obtained through a unified PINN framework, which efficiently handles different types of PDEs, ensuring the system's stability and convergence, especially under complex boundary conditions.

\begin{figure}[htbp] 
    \centering
    \begin{minipage}{0.44\textwidth}
        \centering
        \includegraphics[width=\textwidth]{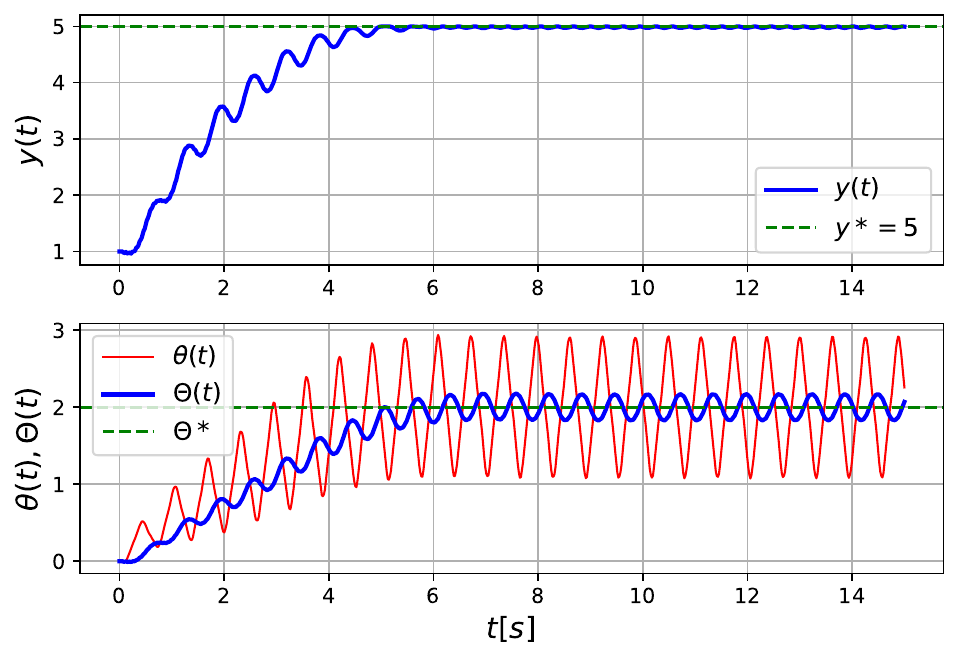}
        \caption{ES results of Distributed PDEs control systems}
        \label{ES1}
    \end{minipage}
    \vspace{0.1cm}  
    \hspace{0.07\textwidth}  
    \begin{minipage}{0.44\textwidth}
        \centering
        \includegraphics[width=\textwidth]{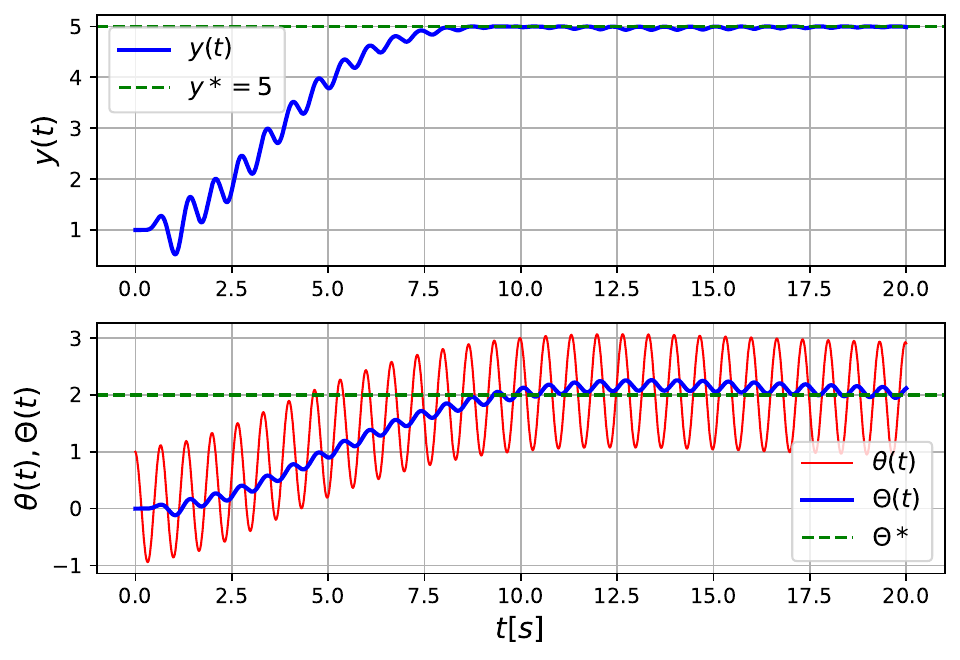}
        \caption{ES results of RAD PDEs control systems}
        \label{ES2}
    \end{minipage}
    \vspace{0.1cm}  
    \centering
    \begin{minipage}[t]{0.44\textwidth}  
        \centering
        \includegraphics[width=\textwidth]{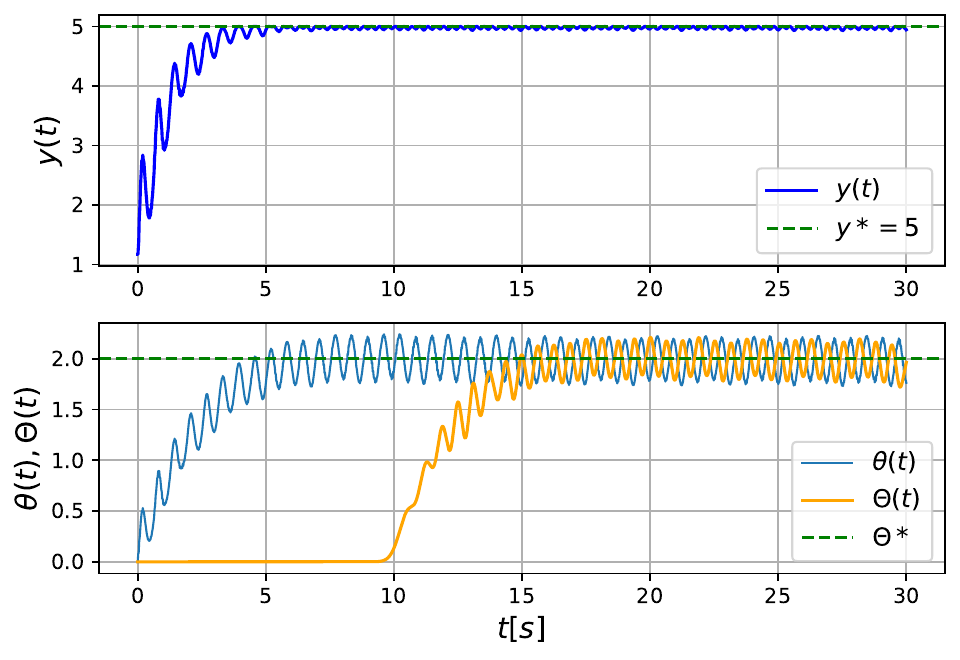}  
        \caption{ES results of Wave PDEs control systems}
        \label{ES3}
    \end{minipage}
    \hspace{0.07\textwidth}  
    \begin{minipage}[t]{0.44\textwidth}  
        \centering
        \includegraphics[width=\textwidth]{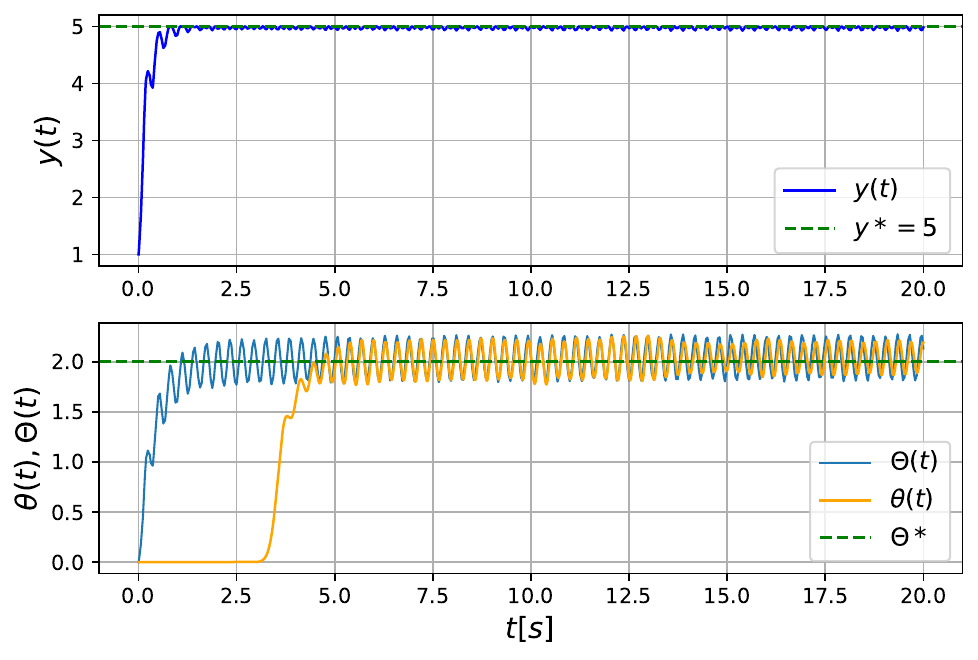}  
        \caption{ES results of Transport PDEs control systems}
        \label{ES4}
    \end{minipage}
\end{figure}

For parabolic PDEs, the static mapping is given in \eqref{05}, and the compensator signal \( G_{\text{av}}(t) \) is designed as shown in \eqref{34}. Figures~\ref{ES1} and~\ref{ES2} present the extremum seeking simulation results for the Distributed PDEs and RAD PDEs systems. In these systems, the input signals and target functions are driven towards the desired values by compensating for the actuator dynamics phase lag. The proposed extremum seeking method combines perturbation signals with an adaptive controller, improving convergence speed and stability. The simulation results demonstrate that, even with phase lag, the controller maintains stability and accelerates convergence to the optimal solution, showcasing its superiority for diffusion PDE systems.

For hyperbolic PDEs, we consider the extremum seeking delay, with the static mapping \( y^\prime(t) \) and compensator signal \( G_\text{av}^\prime(t) \) as in Remark~\ref{remark1}. Figures~\ref{ES3} and~\ref{ES4} show the extremum seeking simulation results for the Wave PDE and Transport Hyperbolic PDE systems, with delays of 10s and 3s, respectively. In both systems, the phase lag is effectively compensated, guiding the input signals and target functions to the desired values. The proposed strategy, combining perturbation signals with an adaptive compensator, compensates for the phase lag in hyperbolic PDEs, improving convergence speed and stability, and enabling rapid convergence to the optimal solution.

\begin{figure}[htbp] 
    \centering
    \begin{minipage}{0.44\textwidth}
        \centering
        \includegraphics[width=\textwidth]{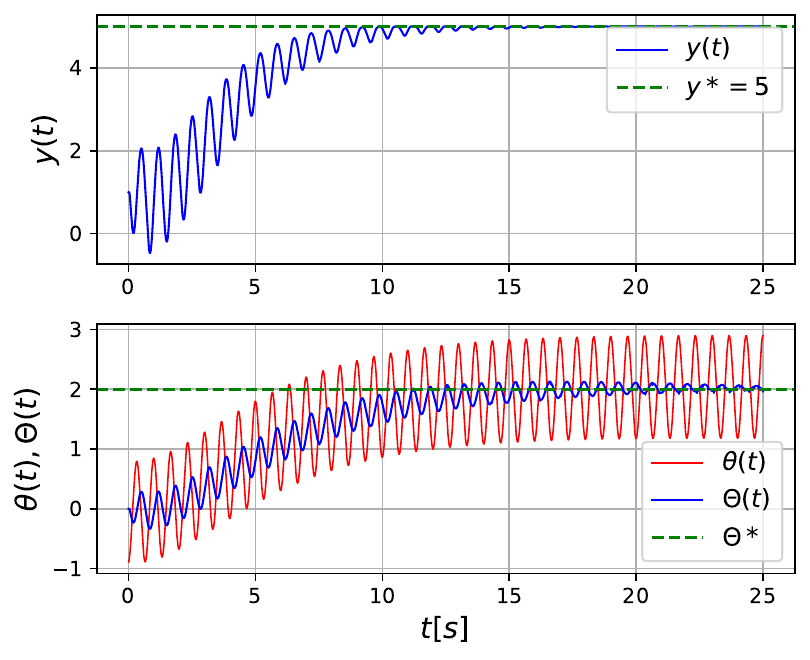}
    \end{minipage}
    \caption{ES results of PDEs with moving boundaries control systems}
    \label{ES5}
\end{figure}

As shown in Fig.~\ref{ES5}, we present the simulation results for PDEs with moving boundaries. We investigate diffusion-type PDEs with moving boundaries, where the boundaries exhibit a decaying trend, significantly affecting the system's dynamic behavior. The simulation results demonstrate that by using the designed compensator, the system can effectively compensate for the changes induced by the moving boundaries, successfully guiding the input signals and target functions towards the desired extremum. Notably, the oscillations of the control inputs and the cost function gradually diminish, and the convergence process becomes more stable, validating the effectiveness of the proposed method.

\section{Conclusions}

This paper proposed a new method for ES of different PDE systems by utilizing PINNs as a unified tool for quickly solving probing/perturbation signals derived from PDEs. By directly embedding physical laws into the neural network training process, we can address additive perturbations signal for various PDE systems within a unified model framework. This eliminates the need to design separate solution strategies for each system. The primary advantage of this method lies in its potential for unity and efficiency. PINNs enable flexibility in adapting to different PDE systems by adjusting inputs and loss functions, solving probing signals in a consistent manner. This approach simplifies the design and implementation of the probing signals, offering notable improvements in computational efficiency. Preliminary results demonstrate that PINNs can effectively approximate the probing signals required for extremum seeking control, providing a promising direction for the application of ES in infinite-dimensional systems. By combining the strengths of deep learning with physical modeling, our approach presents a unified and efficient framework for solving complex PDE systems, with potential extensions to the optimization and control of a broader range of dynamic systems.

Although the framework based on PINNs for solving probing signals derived from PDEs and combined with the ES method has shown promising results, there remain considerable opportunities for further research. As deep learning technologies continue to evolve, a key area for future work will be optimizing the network architecture of PINNs to improve performance in handling complex-valued, nonlinear, and high-dimensional PDE systems (e.g., Schrödinger and Ginzburg–Landau, Kuramoto–Sivashinsky, Korteweg–de Vries, Navier–Stokes equations). Specifically, in the context of extremum seeking optimization control, PDE systems often involve more intricate boundary conditions. Developing more efficient PINN architectures that can quickly and accurately solve the required probing signals will be critical for enhancing the practicality of this method in real-world engineering applications. In addition, future investigation lies in the design and analysis of different control problems with PINNs implementation, as considered in the following references \citep{paper13}, \citep{paper14}, \citep{paper15}, \citep{paper12},
\citep{paper16}, \citep{basa1}, \citep{basa2}, \citep{basa3}.

%
%
%


\begin{thebibliography}{99}

    \bibitem[Adetola \& Guay, 2007]{31}{\scshape Adetola, V. \& Guay, M.}. (2007). Parameter convergence in adaptive extremum-seeking control. {\it Automatica}, {\bf 43}(1), 105--110. 

    \bibitem[Dibo \& Oliveira, 2024]{8}{\scshape Dibo, A. \& Oliveira, T. R.}. (2024). Extremum Seeking Feedback Under Unknown Hessian Signs. {\it IEEE Transactions on Automatic Control}, {\bf 69}, 2383--2390.
    
    \bibitem[Franke \& Schaback, 1998]{1}{\scshape Franke, C. \& Schaback, R.}. (1998). Solving partial differential equations by collocation using radial basis functions. {\it Applied Mathematics and Computation}, {\bf 93}, 73--82.
    
    \bibitem[Galvão et al., 2022]{9}{\scshape Galvão, M. L., Oliveira, T. R. \& Krstić, M.}. (2022). Extremum Seeking for Stefan PDE with Moving Boundary and Delays. {\it IFAC-PapersOnLine}, {\bf 55}, 222--227.
    
    \bibitem[Ghaderi \& Keyanpour, 2022]{2}{\scshape Ghaderi, N. \& Keyanpour, M.}. (2022). Observer-based global output feedback controller design for an unstable wave PDE with nonlinear boundary condition. {\it IMA Journal of Mathematical Control and Information}, {\bf 39}, 295--321.

    \bibitem[Gao et al., 2023]{62}{\scshape Gao, Z., Wang, T., Bai, J., et al.}. (2023). Thermal performance investigation of supercritical methane in minichannel heat sink on flight vehicle actuator under geometry effect of cross section. {\it Numerical Heat Transfer, Part A: Applications}, {\bf 83}(3), 315--330.
    
    \bibitem[Ghaffari et al., 2012]{26}{\scshape Ghaffari, A., Krstić, M. \& Nešić, D.}. (2012). Multivariable Newton-based extremum seeking. {\it Automatica}, {\bf 48}, 1759--1767.
    
    \bibitem[Hale and Lunel, 1990]{30}{\scshape Hale, J. K. \& Lunel, S. M. V.}. (1990). Averaging in infinite dimensions. {\it The Journal of Integral Equations and Applications}, 463--494.
    
    \bibitem[Hudon et al., 2008]{32}{\scshape Hudon, N., Guay, M., Perrier, M., et al.}. (2008). Adaptive extremum-seeking control of convection-reaction distributed reactor with limited actuation. {\it Computers \& Chemical Engineering}, {\bf 32}(12), 2994--3001.

    \bibitem[Hu et al., 2023a]{60}{\scshape Hu, G., Guo, J., Guo, Z., Cieslak, J. \& Henry, D.}. (2023). ADP-Based Intelligent Tracking Algorithm for Reentry Vehicles Subjected to Model and State Uncertainties. {\it IEEE Transactions on Industrial Informatics}, {\bf 19}(4), 6047--6055. 

    \bibitem[Hu et al., 2023b]{61}{\scshape Hu, G., Guo, J., Cieslak, J., Ding, Y., Guo, Z. \& Henry, D.}. (2023). Fault-tolerant control based on adaptive dynamic programming for reentry vehicles subjected to state-dependent actuator fault. {\it Engineering Applications of Artificial Intelligence}, {\bf 123}, 106450. 
    
    \bibitem[Krstić \& Wang, 2000]{10}{\scshape Krstić, M. \& Wang, H.-H.}. (2000). Stability of extremum seeking feedback for general nonlinear dynamic systems. {\it Automatica}, {\bf 36}, 595--601.

    \bibitem[Krstić and Smyshlyaev, 2008]{25}{\scshape Krstić, M. \& Smyshlyaev, A.}. (2008). Boundary control of PDEs: A course on backstepping designs. {\it SIAM}.  

    \bibitem[Krstic, 2009]{24}{\scshape Krstic, M.}. (2009). Compensating actuator and sensor dynamics governed by diffusion PDEs. {\it Systems \& Control Letters}, {\bf 58}, 372--377.

    \bibitem[Lehman, 2002]{35}{\scshape Lehman, B.}. (2002). The influence of delays when averaging slow and fast oscillating systems: overview. {\it IMA Journal of Mathematical Control and Information}, {\bf 19}(1\_and\_2), 201--215.
    
    \bibitem[Li et al., 2017]{3}{\scshape Li, J., Zhai, S., Weng, Z. \& Feng, X.}. (2017). H-adaptive RBF-FD method for the high-dimensional convection-diffusion equation. {\it International Communications in Heat and Mass Transfer}, {\bf 89}, 139--146.
    
    \bibitem[Liu et al., 2020]{4}{\scshape Liu, Y.-Q., Wang, J.-W. \& Sun, C.-Y.}. (2020). A Lyapunov-based design of dynamic feedback compensator for linear parabolic MIMO PDEs. {\it IMA Journal of Mathematical Control and Information}, {\bf 37}, 455--474.
    
    \bibitem[Liu et al., 2024]{11}{\scshape Liu, G., Bai, Y., Zhu, L., Wang, Q. \& Zhang, W.}. (2024). A sequential excitation and simplified ant colony optimization based global extreme seeking control method for performance improvement. {\it Swarm and Evolutionary Computation}, {\bf 86}, 101522.
    
    \bibitem[Lunardi, 2012]{27}{\scshape Lunardi, A.}. (2012). Analytic semigroups and optimal regularity in parabolic problems. {\it Springer Science \& Business Media}.
    
    \bibitem[Manzie \& Krstić, 2009]{12}{\scshape Manzie, C. \& Krstić, M.}. (2009). Extremum Seeking With Stochastic Perturbations. {\it IEEE Transactions on Automatic Control}, {\bf 54}, 580--585.
    
    \bibitem[Meng \& Karniadakis, 2020]{20}{\scshape Meng, X. \& Karniadakis, G. E.}. (2020). A composite neural network that learns from multi-fidelity data: Application to function approximation and inverse PDE problems. {\it Journal of Computational Physics}, {\bf 401}, 109020.
    
    \bibitem[Meng et al., 2023]{21}{\scshape Meng, Z., Qian, Q., Xu, M., Yu, B., Yıldız, A. R. \& Mirjalili, S.}. (2023). PINN-FORM: A new physics-informed neural network for reliability analysis with partial differential equation. {\it Computer Methods in Applied Mechanics and Engineering}, {\bf 414}, 116172.
    
    \bibitem[Nagel, 1989]{29}{\scshape Nagel, R.}. (1989). Towards a “matrix theory” for unbounded operator matrices. {\it Mathematische Zeitschrift}, {\bf 201}, 57--68.

		\bibitem[Nunes et al., 2020]{basa2}{\scshape Nunes, E. V. L., Peixoto, A. J., Oliveira, T. R. \& Hsu, L.}. (2020). Global exact tracking for uncertain multivariable linear systems by output feedback sliding mode control. {\it Proceedings of the 2010 American Control Conference}, {\bf IEEE}, 974--979.

		\bibitem[Oliveira et al., 2017a]{paper14}
{\scshape Oliveira, T.~R., Costa, L.~R., Catunda, J.~M.~Y., Pino, A.~V., Barbosa, W. \& de~Souza, M.~N.}. (2017). Time-scaling based sliding mode control for neuromuscular
  electrical stimulation under uncertain relative degrees. {\it Medical Engineering $\&$ Physics}, {\bf 44} 53--62.
	
		\bibitem[Oliveira et al., 2017c]{paper12}
{\scshape Oliveira, T.~R., Cunha, J.~P. V.~S. \& Hsu, L.}. (2017). Adaptive sliding mode control based on the extended equivalent control concept for disturbances with
  unknown bounds. {\it Advances in Variable Structure Systems and Sliding Mode Control---Theory and Applications. Studies in Systems, Decision and
  Control}, {\bf 115}, 149--163.	
		
		\bibitem[Oliveira et al., 2020]{37}{\scshape Oliveira, T. R., Feiling, J., Koga, S. \& Krstić, M.}. (2020). Multivariable extremum seeking for PDE dynamic systems. {\it IEEE Transactions on Automatic Control}, {\bf 65}(11), 4949--4956. 

    \bibitem[Oliveira et al., 2021]{42}{\scshape Oliveira, T. R., Feiling, J., Koga, S. \& Krstić, M.}. (2021). Extremum seeking for unknown scalar maps in cascade with a class of parabolic partial differential equations. {\it International Journal of Adaptive Control and Signal Processing}, {\bf 35}, 1162--1187.
		
    \bibitem[Oliveira et al., 2019a]{33}{\scshape Oliveira, T. R., Feiling, J., \& Krstic, M.}. (2019). Extremum seeking for maximizing higher derivatives of unknown maps in cascade with reaction-advection-diffusion PDEs. {\it IFAC-PapersOnLine}, {\bf 52}(29), 210--215.		
		
		\bibitem[Oliveira et al., 2017b]{ref}{\scshape Oliveira, T. R., Krstić, M. \& Tsubakino, D.}. (2017). Extremum Seeking for Static Maps With Delays. {\it IEEE Transactions on Automatic Control}, {\bf 62}, 1911--1926.
    
    \bibitem[Oliveira \& Krstić, 2019b]{34}{\scshape Oliveira, T. R. \& Krstić, M.}. (2019). Compensation of wave PDEs in actuator dynamics for extremum seeking feedback. {\it IFAC-PapersOnLine}, {\bf 52}(29), 134--139.

    \bibitem[Oliveira \& Krstić, 2021a]{43}{\scshape Oliveira, T. R. \& Krstić, M.}. (2021). Extremum seeking feedback with wave partial differential equation compensation. {\it Journal of Dynamic Systems, Measurement, and Control}, {\bf 143}(4), 041002.

    \bibitem[Oliveira \& Krstić, 2021b]{13}{\scshape Oliveira, T. R. \& Krstić, M.}. (2021). Extremum seeking boundary control for PDE–PDE cascades. {\it Systems and Control Letters}, {\bf 155}, 105004.

    \bibitem[Oliveira \& Krstić, 2022]{36}{\scshape Oliveira, T. R. \& Krstić, M.}. (2022). Extremum Seeking through Delays and PDEs. {\it SIAM}.
    
    \bibitem[Oliveira \& Krstić, 2023]{14}{\scshape Oliveira, T. R. \& Krstić, M.}. (2023). Extremum seeking for infinite-dimensional systems. {\it Annual Reviews in Control}, {\bf 56}, 100908.
		
		\bibitem[Oliveira et al., 2010]{basa3}{\scshape Oliveira, T. R., Peixoto, A. J., Costa, R. R. \& Hsu, L.}. (2010). Dwell-time and disturbance monitoring for peaking avoidance and performance improvement in high-gain observer based sliding mode control. {\it Dynamics of Continuous, Discrete and Impulsive Systems Series B: Applications and Algorithms}, {\bf 17}, 839--874.		
		
		\bibitem[Oliveira et al., 2015]{paper15}
{\scshape Oliveira, T.~R., Peixoto, A.~J. \& Hsu, L.}. (2015). Global tracking for a class of uncertain nonlinear systems with unknown sign-switching control direction by
  output feedback. {\it International Journal of Control}, {\bf 88}, 1895--1910.
			
		\bibitem[Oliveira et al., 2021]{paper16}
{\scshape Oliveira, T.~R., Rodrigues, V.~H.~P., Krstic, M., \& Basar, T.}. (2021). Nash equilibrium seeking in quadratic noncooperative games under two delayed
  information-sharing schemes. {\it Journal of Optimization Theory and Applications}, {\bf 191}, 700--735.

		\bibitem[Oliveira et al., 2018]{basa1}{\scshape Oliveira, T. R., Rušiti, D., Diagne, M. \& Krstic, M.}. (2018). Gradient extremum seeking with time-varying delays. {\it Proceedings of the 2018 American Control Conference (ACC)}, {\bf IEEE}, 3304--3309.

    \bibitem[Pazy, 2012]{28}{\scshape Pazy, A.}. (2012). Semigroups of linear operators and applications to partial differential equations. {\it Springer Science \& Business Media}.
		
		\bibitem[Pinto et al., 2019]{paper13}
{\scshape Pinto, H.~L. C.~P., Oliveira, T.~R. \& Hsu, L.}. (2019). Sliding mode observer for fault reconstruction of time-delay and sampled-output systems---a time shift
  approach. {\it Automatica}, {\bf 106}, 390--400.

    \bibitem[Radenković et al., 2018]{15}{\scshape Radenković, M. S., Stanković, M. S. \& Stanković, S. S.}. (2018). On Stochastic Extremum Seeking via Adaptive Perturbation–Demodulation Loop. {\it Journal of Optimization Theory and Applications}, {\bf 179}, 1008--1024.

    \bibitem[Raissi et al., 2019]{22}{\scshape Raissi, M., Perdikaris, P. \& Karniadakis, G. E.}. (2019). Physics-informed neural networks: A deep learning framework for solving forward and inverse problems involving nonlinear partial differential equations. {\it Journal of Computational Physics}, {\bf 378}, 686--707.

    \bibitem[Rušiti et al., 2019]{38}{\scshape Rušiti, D., Evangelisti, G., Oliveira, T. R., Gerdts, M. \& Krstić, M.}. (2019). Stochastic extremum seeking for dynamic maps with delays. {\it IEEE Control Systems Letters}, {\bf 3}(1), 61--66. 

    \bibitem[Rušiti et al., 2021a]{39}{\scshape Rušiti, D., Oliveira, T. R., Krstić, M. \& Gerdts, M.}. (2021). Newton-based extremum seeking of higher-derivative maps with time-varying delays. {\it International Journal of Adaptive Control and Signal Processing}, {\bf 35}, 1202--1216.

    \bibitem[Rušiti et al., 2021b]{40}{\scshape Rušiti, D., Oliveira, T. R., Krstić, M. \& Gerdts, M.}. (2021). Robustness to delay mismatch in extremum seeking. {\it European Journal of Control}, {\bf 62}, 75--83. 
    
    \bibitem[Ruthotto \& Haber, 2020]{23}{\scshape Ruthotto, L. \& Haber, E.}. (2020). Deep Neural Networks Motivated by Partial Differential Equations. {\it Journal of Mathematical Imaging and Vision}, {\bf 62}, 352--364.
    
    \bibitem[Silva et al., 2024a]{5}{\scshape Silva, P. C. S., Pellanda, P. C., Oliveira, T. R., Andrade, G. A. D. \& Krstic, M.}. (2024). Extremum Seeking for a Class of Wave Partial Differential Equations With Kelvin-Voigt Damping. {\it IEEE Control Systems Letters}, {\bf 8}, 43--48.
    
    \bibitem[Silva et al., 2024b]{16}{\scshape Silva, P. C. S., Pellanda, P. C. \& Oliveira, T. R.}. (2024). Stochastic Multivariable Extremum Seeking Control Considering Input and Output Delays. {\it Journal of Control, Automation and Electrical Systems}, {\bf 35}, 986--998.

    \bibitem[Tadmor, 2012]{6}{\scshape Tadmor, E.}. (2012). A review of numerical methods for nonlinear partial differential equations. {\it Bulletin of the American Mathematical Society}, {\bf 49}, 507--554.

    \bibitem[Tan et al., 2010]{17}{\scshape Tan, Y., Moase, W. H., Manzie, C., Nešić, D. \& Mareels, I. M. Y.}. (2010). Extremum seeking from 1922 to 2010. {\it Proceedings of the 29th Chinese Control Conference}, 29--31 July 2010, 14--26.

    \bibitem[Tsubakino et al., 2023]{41}{\scshape Tsubakino, D., Oliveira, T. R. \& Krstić, M.}. (2023). Extremum seeking for distributed delays. {\it Automatica}, {\bf 153}, 111044.

    \bibitem[Yuan et al., 2024]{63}{\scshape Yuan, R., Guo, Z., Cao, S., Henry, D., Cieslak, J., Oliveira, T. R. \& Guo, J.}. (2024). Physics-Informed Neural Networks-based Uncertainty Identification and Control for Closed-Loop Attitude Dynamics of Reentry Vehicles. {\it Proceedings of the 2024 IEEE 18th International Conference on Control \& Automation (ICCA)}, 629--634.
    
    \bibitem[Zhai et al., 2015]{7}{\scshape Zhai, S., Qian, L., Gui, D. \& Feng, X.}. (2015). A block-centered characteristic finite difference method for convection-dominated diffusion equation. {\it International Communications in Heat and Mass Transfer}, {\bf 61}, 1--7.
				
		 
\end{thebibliography}
\end{document}